\documentclass[11pt]{amsart}
\usepackage{verbatim, latexsym, amssymb, amsmath }
\def\R{\mathbb R}
\def\N{\mathbb N}

\def\h2{\mathrm{area}}
\def\v{\mathrm{volume}}

\def\area{\left\lvert\Sigma\right\rvert}
\def\a_t{\left\lvert\Sigma_t\right\rvert}
\def\b_i{\left\lvert\Sigma_i\right\rvert}
\def\d{\mathrm{div}}
\newcommand{\fint}{\mathop{\int\makebox(-15,2){\rule[4pt]{.9em}{0.6pt}}\kern-4pt}\nolimits}
\def\tf {\mathring{A}}
\def\n{\lvert\mathring{A}\rvert}

\def\dr{\lvert\partial_r^{\top}\rvert}
\def\dw_i{\lvert\widehat{\nabla} w_t\rvert}
\def\ge{\lvert{\nabla}_0 w\rvert}

\def\as {asymptotically }

\newtheorem{thm}{Theorem}[section]
\newtheorem{lemm}[thm]{Lemma}
\newtheorem{cor}[thm]{Corollary}
\newtheorem{prop}[thm]{Proposition}
\theoremstyle{remark}
\newtheorem{rmk}[thm]{Remark}
\theoremstyle{definition}
\newtheorem{defi}[thm]{Definition}

\title{Existence and Uniqueness of constant mean curvature foliation of asymptotically hyperbolic $3$-manifolds}

\author{Andr\'e Neves}
\email{aneves@math.princeton.edu}
\address{Fine Hall, Princeton University,
Princeton, NJ 08544, USA}

\author{Gang Tian}
\email{tian@math.princeton.edu}
\address{Fine Hall, Princeton University,
Princeton, NJ 08544, USA}

\begin{document}

\maketitle \markboth{Existence and uniqueness of constant mean curvature foliations} { Andr\'e Neves and Gang
Tian} \maketitle

\begin{abstract}
We prove existence and uniqueness of foliations by stable spheres  with constant mean curvature for
$3$-manifolds which are asymptotic to  Anti-de Sitter-Schwarzschild metrics with positive mass. These metrics
arise naturally as spacelike timeslices for solutions of the Einstein equation with a negative cosmological
constant.
\end{abstract}

\section{Introduction}

In 1996, Huisken and Yau \cite{Huisken3} showed that metrics that arise as  spacelike timeslices of solutions to
the Einstein vacuum equation have a well defined center of mass. More precisely, they showed that metrics which
are asymptotic to Schwarzschild metrics with positive mass outside a compact set admit, under some technical assumptions, a unique foliation by
stable spheres with constant mean curvature.  We prove an analogous result for a large class of metrics that
arise naturally as spacelike timeslices for solutions to the Einstein equation with negative cosmological
constant. More precisely, we show that metrics which are, outside a compact set, asymptotic to Anti--de
Sitter--Schwarzschild metrics with positive mass admit, under some technical assumptions, a unique foliation by stable spheres with constant mean
curvature and thus, they have a well defined center of mass.

Most of the work in this paper is devoted to prove uniqueness of foliations by stable spheres with constant mean
curvature for \as Anti--de Sitter--Schwarzschild metrics with positive mass. Note that the uniqueness of such
foliations fails for the hyperbolic $3$-space and so the positivity of the mass needs to be used. The central
idea in the paper consists in combining the positivity of the mass with the Kazdan-Warner obstructions
\cite{Kazdan} in order to prove uniqueness. The relation between positive mass and the study of foliations near infinity was  observed by Christodoulou and Yau in \cite{CK}.

The existence of such foliations in the \as flat setting was proven by Huisken and Yau \cite{Huisken3} using a
modified mean curvature flow and by Ye \cite{Ye} using a perturbation method. For metrics asymptotic to Anti--de
Sitter--Schwarzschild metrics, the existence of such foliations was proven by Rigger \cite{Rigger} using the
mean curvature flow approach. The
 arguments in \cite{Ye} can be adapted in a straightforward way in order to obtain the existence result for
metrics asymptotic to Anti--de Sitter--Schwarzschild metrics. We include its proof in the last section for the
sake completeness.

Before stating the main result we need to introduce some notation.
 Denote by $g_0$ the standard round
metric on $S^2$. The Anti--de Sitter--Schwarzschild manifold with mass $m>0$ is defined to be
$(s_0,\infty)\times S^2$ with the metric
$$g_m=(1+s^2-m/s)^{-1}ds^2+s^2g_0,$$
where $s_0$ is the zero of $1+s^2-m/s$. Note that when $m=0$  the Anti--de Sitter--Schwarzschild metric becomes
the hyperbolic metric.  After a change of coordinates, the metric can be written in the form
$$g_m=dr^2+(\sinh^2r+m/(3\sinh r)+O(\exp(-3r)))g_0.$$

Our result will apply to metrics which are, outside a compact set, lower order perturbations of Anti--de
Sitter--Schwarzschild metrics.
\begin{defi}\label{dd}
$(M,g)$ is an \as Anti--de Sitter--Schwarzschild manifold with mass $m$ if, for some compact set $K$, $M-K$ is
diffeomorphic to $\R^3$ minus a ball and, with respect to this diffeomorphism, the metric can be written (in
spherical coordinates) as
$$g=dr^2+(\sinh^2r+m/(3\sinh r))g_0+Q,$$
where
$$\lvert Q\rvert+\lvert \nabla Q\rvert +\lvert \nabla^2 Q\rvert=O(\exp(-5r)).$$
\end{defi}
\begin{rmk}
A direct computation shows that the mass $m$ coincides with the mass defined in \cite{herzlich} and
\cite{wang} for \as hyperbolic metrics.
\end{rmk}
With respect to the coordinates specified in the Definition \ref{dd}, $M-K$ becomes equipped with a radial
function $r$. Hence, given a foliation $(\Sigma_t)_{t\geq 0}$ of $M$, we can define the lower radius and the
upper radius to be
$$\underline{r}_t=\sup\{r(x)\,|\, x\in \Sigma_t\}\quad\mbox{and}\quad\overline{r}_t=\inf\{r(x)\,|\, x\in \Sigma_t\}$$
respectively. Given  a family of functions $f_t$ defined on $\Sigma_t$, we use
$$f_t=O(\exp(-n\underline{r}_t))\quad\mbox{and}\quad f_t=o(\exp(-n\underline{r}_t))$$ to denote that
$$\limsup_{t\to\infty}(\vert f_t \rvert \exp(n\underline{r}_t))<
\infty\quad\mbox{and}\quad\limsup_{t\to\infty}(\vert f_t \rvert \exp(n\underline{r}_t))=0$$ respectively.

A surface $\Sigma$ with constant mean curvature is said to be stable if volume preserving variations do not
decrease its area. A standard computation shows that stability is equivalent to the second variation operator
$$ Lu\equiv-\Delta f -\left(\lvert A\rvert^2+R(\nu,\nu) \right)f$$
having only nonnegative eigenvalues when restricted to functions with zero mean value, i.e.,
     $$\int_{\Sigma_t}\left(\lvert A\rvert^2+R(\nu,\nu) \right)f^2d\mu\leq \int_{\Sigma_t}\lvert \nabla f \rvert^2 d\mu$$
for all functions $f$ with $\int_{\Sigma_t}fd\mu=0.$

The purpose of this paper is to show the following  result.
\begin{thm}\label{main1}
    Asymptotically Anti--de Sitter--Schwarzschild manifolds with positive mass admit a unique foliation
    $(\Sigma_t)_{t\geq 0}$ by stable spheres with constant mean curvature such that
    \begin{equation}\label{condi}
            \lim_{t\to\infty}(\overline{r}_t-6/5\underline{r}_t)=-\infty.
    \end{equation}
\end{thm}
\begin{rmk}
A condition similar to  \eqref{condi}  was also assumed in \cite{Huisken3} for metrics that are asymptotic
Schwarzschild metrics. It is an interesting question if one can strengthen Theorem \ref{main1} by weakening
condition \eqref{condi} to the condition that $\underline{r}_t$ is sufficiently large. Its analogous version for
metrics asymptotic to Schwarzschild metrics was solved in \cite{qingtian}.
\end{rmk}

A contradiction argument implies the following corollary.
\begin{cor}
There are positive constants $C_0$ and $R_0$ depending only on $g$ so that any stable sphere $\Sigma$ with
constant mean curvature satisfying
$$\lvert\Sigma\rvert\geq C_0\quad\mbox{and}\quad\overline{r}-6/5\underline{r}\leq-R_0$$
is unique.
\end{cor}
\begin{proof}
    Assume there are two distinct sequences $(\Sigma^1_i)_{i\in\N}$ and $(\Sigma^2_i)_{i\in\N}$ of stable
    spheres with constant mean curvature such that
    \begin{align*}
    &\lim_{i\to\infty}\lvert\Sigma^1_i\rvert=\lim_{i\to\infty}\lvert\Sigma^2_i\rvert=\infty,\\
    &\lim_{i\to\infty}\overline{r}^1_i-6/5\underline{r}^1_i=\lim_{i\to\infty}\overline{r}^2_i-6/5
    \underline{r}^2_i=-\infty,
    \end{align*}
    and $H(\Sigma^1_i)=H(\Sigma^2_i)$ for every integer $i$. The proof of Theorem \ref{main1} also applies to the
    foliations
    $$\Sigma^j_t\equiv\Sigma^j_i\quad\mbox{if }i\leq t<i+1,$$
    where $j=1,2,$ and so we obtain a contradiction.
\end{proof}

{\bf Acknowledgement} The first author would like to express his gratitude to Alice Chang for many useful
discussions.

\section{Outline of the proof}
We outline the proof of Theorem \ref{main1} in order to emphasize the main ideas over the technical aspects of
the paper. Let $(\Sigma_t)_{t\geq 0}$ denote a foliation by stable spheres with constant mean curvature. For the
sake of simplicity, we assume in this discussion that $\overline{r}_t-\underline{r}_t$ is uniformly bounded. We
use the notation
$$\fint_{\Sigma_t}f\,d\mu \equiv |\Sigma_t|^{-1}\int_{\Sigma_t}f\,d\mu,$$
where $|\Sigma_t|$ stands for the surface area of $\Sigma_t$.

 Section \ref{prelim} is devoted to auxiliary computations. In
Section \ref{intestimate} we follow the same argumentation done by Huisken and Yau in \cite{Huisken3} and use
the stability assumption in order to derive the following estimate for the mean curvature of $\Sigma_t$ (Lemma
\ref{h})
\begin{equation}\label{mozart}
    H^2=4+16\pi/\lvert\Sigma_t\rvert+\fint_{\Sigma_t}O(\exp(-3r))d\mu
\end{equation}
 and the following integral estimate for the
trace free part of the second fundamental form (Proposition \ref{integral})
$$\fint_{\Sigma_t}\n^2d\mu\leq O(\exp(-4\underline{r}_t)).$$

In Section \ref{Intrin} we study the intrinsic geometry of $\Sigma_t$. More precisely, we show that after
pulling back by a suitable diffeomorphism from $\Sigma_t$ to $S^2$, the metric
$$\hat{g}_t\equiv 4\pi/|\Sigma_t|g$$
can be written as
$$\exp(2\beta_t)g_0,$$
where $g_0$ denotes the standard round metric on $S^2$ and $$\beta_t=O(\exp(-\underline{r}_t)).$$ This result
implies that $\hat{g}_t$ is very ``close'' to being a round metric. The proof of this result (Theorem
\ref{intrinsic}) has two steps.

The first step consists in deriving a pointwise estimate for $\n$ from the integral estimate (Proposition
\ref{a}). In order to do so, we have to exploit the fact that the hyperbolic metric is conformal to the
Euclidean metric on the unit ball and so the same is true, up to a term of low order, for the metric $g$.
Therefore, denoting by $d\sigma$ the surface measure induced by the Euclidean metric on $\Sigma_t$, we have by
conformal invariance that
$$ \lim_{t\to\infty}\int_{\Sigma_t}\overline{\n}^2 d\sigma=\lim_{t\to\infty}\int_{\Sigma_t}{\n}^2\,d\mu=0,$$
where the quantities measured with respect to the Euclidean metric are denoted with a bar. Because the Euclidean
area of $\Sigma_t$ converges to $4\pi$ (Proposition \ref{estimates}), the identity above and Gauss-Bonnet
Theorem imply that the Euclidean mean curvature $\bar{H}$ has no concentration points. We can then use
Michael-Simon Sobolev inequality and the equation satisfied by $\n^2$ (with respect to the Euclidean metric) in
order to apply the standard Moser iteration procedure and conclude that, up to lower order terms,
\begin{equation}\label{schubert}
\sup_{\Sigma_t}\n^2\leq C\int_{\Sigma_t}\n^2\,d\sigma=C\fint_{\Sigma_t}\n^2d\mu\leq O(\exp(-4\underline{r}_t)).
\end{equation}

The second step consists in using Gauss equation (Lemma \ref{ricci}) which, combined with estimates
\eqref{mozart} and \eqref{schubert}, implies that the Gaussian curvature of $\Sigma_t$ with respect to $\hat
g_t$ satisfies
        $$\hat K_t=(4\pi)^{-1}|\Sigma_t|((H^2-4)/4  -\lvert\mathring{A}\rvert^2/2)+O(\exp(-r))=1+O(\exp(-r)).$$
The desired result follows from this estimate.

Section \ref{Unique} contains the main estimate that makes uniqueness possible. Set
$$w_t(x)\equiv r(x)-\hat{r}_t,\quad\mbox{ where }\quad \a_t=4\pi\sinh^2\hat{r}_t.$$ We want to show that $w_t$
converges uniformly to zero (Theorem \ref{unique}). With respect to the standard round metric $g_0$, the
functions $w_t$ satisfy the equation (see \eqref{e2})
$${\Delta}_0 w_t = \exp(-2w_t)- 1+P,$$
where $$\int_{S^2}|P|d\mu_0=O(\exp(-r)).$$ Because we are assuming that $\overline{r}_t-\underline{r}_t$ is
uniformly bounded, we have that $w_t$ is bounded in $W^{1,2}$ with respect to $g_0$ and thus, we can take a
sequence converging weakly to $w_0$ that satisfies
\begin{equation}\label{jazz}
{\Delta}_0 w_0 = \exp(-2w_0)- 1.
\end{equation}
 It is well known that this equation can have many solutions. Therefore, we need
to use the fact that the mass is nonzero in order to show that $w_0=0$. This is achieved through the
Kazdan-Warner identity \cite{Kazdan}. Because the metric $\hat g_t$ is ``close'' to being the round metric, this
identity implies that, for each of the standard coordinate functions $x_1,x_2,x_3$ on $S^2$ (see
\eqref{brubaker}),
$$\int_{S^2}x_i\hat{K}_t d{\mu_0}=0,\quad\mbox{ for } i=1,2,3.$$
On the other hand, a careful expansion of the terms involved in the Gauss equation shows that the Gaussian
curvature of $\hat{g}_t$ is such that
\begin{multline*}
\a_t^{1/2}\hat{K}_t= (4\pi)^{-1}\a_t^{3/2}(H^2-4)/4\\
+m(4\pi)^{-1}\a_t^{3/2} /{\sinh^3 r}+O(\exp(-r)).
\end{multline*}
Hence, because the mean curvature is constant, we obtain from the Kazdan-Warner identity that (Proposition
\ref{kw})
$$
 m\int_{S^2} x_i \exp(-3 w_t) d{\mu}_0=O(\exp(-\underline{r}_t)),\quad\mbox{ for } i=1,2,3.
$$
Therefore,
$$
m\int_{S^2} x_i \exp(-3 w_0) d{\mu}_0=0,\quad\mbox{ for } i=1,2,3.
$$
Recalling that $w_0$ is a solution to equation \eqref{jazz}, the above identity implies that $w_0=0$. Standard
techniques  can then be used to show that $w_t$ converges to zero uniformly.

In Section \ref{strong} we improve the rate of convergence of $w_t$ to zero (Theorem \ref{close}). In order to
do so, we redefine $w_t$ to be
$$w_t(x)\equiv r(x)-\tilde{r}_t,$$
where $\tilde{r}_t$ is such that
 $$H=2\cosh \tilde{r}_t/\sinh{\tilde{r}_t}-m /\sinh^3{\tilde{r}_t}+o(\exp(-4\tilde{r_t})).$$
Then, the equation satisfied by $w_t$ improves to become
$${\Delta}_0 w_t = \exp(-2w_t)- 1+o(\exp(-2r)).$$ Using the orthogonality condition given by the Kazdan-Warner
obstructions, standard elliptic estimates  show that
 $$\lvert w_t \rvert_{C^{2,\alpha}}=o(\exp(-\underline{r}_t)).$$

In Section \ref{uni}, we use the strong approximation of $\Sigma_t$ to a coordinate sphere in order to prove
uniqueness of foliations by stable spheres with constant mean curvature. The main reason for this to work is
that, with respect to the round metric on $S^2$, the linearization of the mean curvature on a coordinate sphere
$\{\lvert x\rvert=r\}$ is the operator
$$L(f)=\Delta_0 f+(2-3m/\sinh r)f,$$
which is invertible is $m$ is not zero. Finally, in Section \ref{exi} we adapt the arguments used in \cite{Ye}
and we show that, for all $r$ sufficiently large, we can find a stable sphere with constant mean curvature which
is a perturbation of $\{\lvert x\rvert=r\}$.

\section{Preliminaries}\label{prelim}

In this section we compute the relevant formulas needed throughout this paper. Before doing so, we need to
introduce some notation. $\Sigma$ is assumed to be a stable sphere with constant mean curvature. The radial
vector is denoted by $\partial_r$ and $\partial_r^{\top}$ stands for the tangential projection of $\partial_r$
on $T\Sigma$, which has length denoted by $\dr$. Finally, $\nu$ denotes the exterior unit normal to $\Sigma$.

We start by computing the asymptotic expansion of some geometric quantities depending on $g_m$, the Anti de
Sitter- Schwarzschild metric with mass $m$. Let $\{e_1,e_2\}$ denote a $g_m$-orthonormal basis for the
coordinate spheres $\{\lvert x\rvert=r\}.$
\begin{lemm}\label{ads}
$ $
    \begin{enumerate}
        \item[(i)] The mean curvature $H_m(r)$ of $\{\lvert x \rvert=r\}$ is such that
            $$H_m(r)=2\cosh r/\sinh r-m/\sinh^3 r+O(\exp(-5r)).$$
        \item[(ii)] The scalar curvature is -6.
        \item[(iii)] The Ricci curvature is such that
            \begin{align*}
                R_m(\partial_r,\partial_r)&=-2-{m}/{\sinh^3 r}+O(\exp(-5r))\\
                R_m(e_1,e_1)&=R_m(e_2,e_2)=-2+{m}/({2\sinh^3 r})+O(\exp(-5r))\\
                R_m(e_1,e_2)&=0.
            \end{align*}
        \item[(iv)]
            The derivatives of the Ricci curvature with respect $g_m$ are such that
                $$ \nabla^m_{e_1}R_m(\partial_r,e_1) =\nabla^m_{e_2}R_m(\partial_r,e_2).$$
        \end{enumerate}
\end{lemm}
\begin{proof}
If we write $g_m$ as $dr^2+\psi^2(r)g_0$, then $$H_m(r)=2\psi'(r)/\psi(r)$$ and so the first formula follows.
The second formula is just direct computation. The rotational symmetry of the metric implies that its Gaussian
curvature is $\psi^{-2}(r)$ and that the second fundamental form of $\{\lvert x \rvert=r\}$ is trace free.
Hence, we have from Gauss equation that
 $$R_m(\partial_r,\partial_r)=-3-\psi^{-2}(r)+H^2/4$$ and so the first identity in (iii) follows. The
 other two identities in (iii) are a consequence of (ii) and rotational symmetry. The last identity follows from
 the same type of arguments.
\end{proof}

The next lemma  relates the mass $m$ of a metric $g$ with the Gaussian curvature $K$ of a surface $\Sigma$.
\begin{lemm}\label{ricci} The Gaussian curvature of $\Sigma$ satisfies
    \begin{multline*}
        K=(H^2-4)/4+m/\sinh^3 r\\-3m \dr^2/{(2\sinh^3 r})-  \lvert\mathring{A}\rvert^2/2+O(\exp(-5r)).
    \end{multline*}
\end{lemm}
\begin{proof}
    Because $g$ is $C^2$-perturbation of order $O(\exp(5r))$ of $g_m$, we have
        $$R=-6+O(\exp(-5r))$$
    and
\begin{align*}
    2+R(\nu,\nu) & =2+R_m(\nu,\nu)+O(\exp(-5r))\\
    &=-{m}/{\sinh^3 r}+3m \dr^2/{(2\sinh^3 r})+O(\exp(-5r)).
\end{align*}
 The result follows from  Gauss equation
$$K=R/2-R(\nu,\nu)+H^2/4-\lvert \mathring{A}\rvert^2 /2.$$
\end{proof}

Next, we derive the equation satisfied by the trace-free part of the second fundamental form ${\mathring{A}}$.
\begin{lemm}\label{A}
    The Laplacian of $|\tf|^2$ satisfies
    \begin{multline*}
         \Delta ({\lvert\mathring{A}\rvert^2}/{2}) = \left(\frac{H^2-4}{2}-
         \lvert \tf \rvert ^2 \right)\lvert \tf \rvert ^2+\lvert\nabla\tf\rvert^2 \\
         + \left( \lvert\tf\rvert^2 +\lvert H \rvert\lvert \partial _r^{\top} \rvert^2\n\right)
         O(\exp(-3r))+\n O(\exp(-5r)).
    \end{multline*}
\end{lemm}

\begin{proof}

    We assume normal coordinates $x=\{x^i\}_{i=1,2}$ around a point $p$ in the constant mean curvature surface
    $\Sigma$. The tangent vectors are denoted by $\{\partial_1,\partial_2\}$, the normal vector by $\nu$, and the
    Einstein summation convention for the sum of repeated indices is used.

    Simons' identity for the Laplacian of the second fundamental form $A$ (see for instance
    \cite{Huisken1}) implies that
    \begin{multline*}
        \Delta (|\mathring A|^2/2) = \lvert\nabla\tf\rvert^2+H\mbox{Tr}(\tf^3)+H^2\n^2/2-\n^4\\
        +H\tf_{ij}{R}_{\nu i \nu j}-{R}_{\nu\nu}\n^2
        +2{R}_{kikm}\tf_{mj}\tf_{ij}\\
        +2{R}_{kijm}\tf_{km}\tf_{ij}+ \tf_{ij}({\nabla}_k {R}_{\nu j i k}+ {\nabla}_i {R}_{\nu k j k}).
    \end{multline*}
    From
    \begin{equation*}
        {R}_{stuv}=-(\delta_{su}\delta_{tv}-\delta_{sv}\delta_{tu})+O(\exp({-3r}))
    \end{equation*}
    if follows that
        $$2{R}_{kikm}\tf_{mj}\tf_{ij}+2{R}_{kijm}\tf_{km}\tf_{ij}=-4\n^2+\n^2O(\exp(-3r)).$$
    Hence, $\mbox{Tr}(\tf^3)=0$ implies that
    \begin{multline*}
        \Delta({\lvert\mathring{A}\rvert^2}/{2})  = \lvert\nabla\tf\rvert^2+\left(\frac{H^2-4}{2}
        -\lvert \tf \rvert ^2 \right)\lvert \tf \rvert ^2+ H{R}_{\nu i \nu j}\tf_{ij}\\
         + ({\nabla}_k {R}_{\nu j i k}+{\nabla}_i {R}_{\nu k j k})\tf_{ij}+ \lvert\tf\rvert^2 O(\exp(-3r)).
    \end{multline*}
    On the other hand, if $\{v_1,v_2\}$ is an eigenbasis for $\tf$, we obtain from Lemma \ref{ads} that
    \begin{align*}
        {R}_{\nu i \nu j}\tf_{ij}& =\frac{\n}{\sqrt 2}({R}(v_1,v_1)-{R}(v_2,v_2))\\
        & = \frac{\n}{\sqrt 2}({R}_m(v_1,v_1)-{R}_m(v_2,v_2))+\n O(\exp(-5r))\\
        &=\lvert \partial _r^{\top} \rvert^2\n O(\exp(-3r))+\n O(\exp(-5r))
    \end{align*}
    and
    \begin{align*}
        ({\nabla}_k  {R}_{\nu j i k}+ {\nabla}_i  {R}_{\nu k j k})\tf_{ij} & = \sqrt{2}\n( {\nabla}_{v_1} {R}(\nu,v_1)-
        {\nabla}_{v_2} {R}(\nu,v_2))\\
        & = \lvert \partial _r^{\top}\rvert^2\n O(\exp(-3r))+\n O(\exp(-5r)).
    \end{align*}
    Thus, the desired result follows.
\end{proof}

Finally, we derive the equation satisfied by the Laplacian of $r$ on $\Sigma$. Throughout the rest of this
paper, we will progressively estimate and explore all of its terms.

\begin{prop}\label{laplace}
The Laplacian of $r$ on $\Sigma$ satisfies
\begin{multline*}
    \Delta r = (4-2\lvert\partial _r^{\top} \rvert^2 )\exp(-2r)+2-H\\ +(H-2)(1-\langle \partial _r, \nu \rangle)+(1-\langle \partial _r, \nu \rangle)^2
        + O(\exp(-3r))
\end{multline*}
or, being more detailed,
\begin{multline*}
    \Delta r = 2{\cosh r}/{\sinh r}- {m}/{\sinh^3 r}-H \\+(H-2)(1-\langle \partial _r, \nu \rangle)
    +(1-\langle \partial _r, \nu \rangle)^2
        -2\dr^2\exp(-2r)\\+\dr^2O(\exp{(-3r)})+O(\exp(-5r)).
\end{multline*}
\end{prop}
\begin{proof}
    It suffices to prove the second identity. Using again normal coordinates $\{x_1,x_2\}$, we have from Lemma \ref{ads}
            \begin{align*}
                \d_{\Sigma}\partial_r & =  \langle\nabla^m_{i}\partial_r,\partial_i\rangle+O(\exp(-5r))\\
                 & =  H^m(r)-H^m(r)\dr^2/2+O(\exp(-5r))\\
                 & =  (2-\dr^2){\cosh r}/{\sinh r}-{m}/{\sinh^3 r}\\&\quad+\dr^2O(\exp{(-3r)})+ O(\exp(-5r)).
            \end{align*}
    On the other hand,
        $$\d_{\Sigma}\partial_r=\Delta r+ \langle \partial_r , \nu\rangle H $$  and the result follows from the easily
        checked identity
        \begin{multline*}
            \dr^2{\cosh r}/{\sinh r}+\langle \partial_r , \nu\rangle H=H -(H-2)(1-\langle \partial _r, \nu \rangle)\\
            -(1-\langle \partial _r, \nu \rangle)^2+2\lvert\partial _r^{\top}
            \rvert^2\exp(-2r)+\dr^2O(\exp{(-4r)}).
        \end{multline*}
\end{proof}

\section{Integral estimates}\label{intestimate}

We use the stability condition in the same spirit as in \cite{Huisken3}  in order to estimate the mean curvature
$H$ and the $L^2$ norm of $\n$. In the next section,  we combine these integral estimates with Lemma \ref{A} in
order to obtain pointwise estimates for $\n$. From this section on, $$(\Sigma_t)_{t\geq 0}$$ denotes a foliation
by stable spheres with constant mean curvature satisfying condition \eqref{condi}.
We omit the index $t$ in the notation whenever it becomes obvious that we are referring to quantities depending
on $\Sigma_t$.

\begin{lemm}\label{h}
    For each $t$, we have that
        $$H^2=4+16\pi/\lvert\Sigma_t\rvert+\fint_{\Sigma_t}O(\exp(-3r))d\mu$$
    or, equivalently,
        $$H=2+4\pi/\lvert\Sigma_t\rvert+\fint_{\Sigma_t}O(\exp(-3r))d\mu. $$
\end{lemm}
\begin{proof}
    The stability condition implies that, after a clever choice of test function (see \cite[Proposition 5.3]{Huisken3}),
        $$8\pi\geq\int_{\Sigma_t}\lvert A\rvert^2+R(\nu,\nu)d\mu=\int_{\Sigma_t}\n^2+\left(H^2-4\right)/2+R(\nu,\nu)+2\,d\mu$$
    and so, because
        $$R(\nu,\nu)=-2+O(\exp(-3r)),$$
     we have
        $$H^2\leq 4+16\pi/\lvert\Sigma_t\rvert+\fint_{\Sigma_t}O(\exp(-3r))d\mu.$$
    On the other hand, from  Lemma \ref{ricci} and Gauss-Bonnet Theorem  we obtain that
        $$H^2=4+16\pi/\area+2\fint_{\Sigma_t}\n^2d\mu+\fint_{\Sigma_t}O(\exp(-3r))d\mu,$$
    and thus the result follows.
\end{proof}

This lemma combined with the equation for $\Delta r$ gives us these first estimates regarding $\Sigma_t$.

\begin{prop}\label{estimates}
 The following identities hold:
    \begin{enumerate}
        \item[(i)]
            $$\int_{\Sigma_t}\exp{(-2r)}d\mu=\pi+O(\exp(-\underline{r}_t)).$$
        In particular, there is a constant $C$ so that $$C^{-1}\exp(2\underline{r}_t)\leq \a_t\leq C\exp(2\overline{r}_t).$$
        \item[(ii)]$$\int_{\Sigma_t}(1-\langle\partial_r,\nu\rangle)^2d\mu=O(\exp(-\underline{r}_t)).$$
        \item[(iii)]For every positive integer $j$
                    $$\int_{\Sigma_t}j\dr^2\exp(-jr)d\mu = O(\exp(-j\underline{r}_t)).$$
    \end{enumerate}
\end{prop}
\begin{proof}
    Integrating the first identity in Proposition \ref{laplace} and using Lemma \ref{h} we obtain
        \begin{multline}\label{e1}
             4\pi+ \int_{\Sigma_t}O(\exp(-3r))d\mu  =
            \int_{\Sigma_t}(4-2\lvert\partial _r^{\top} \rvert^2 )\exp(-2r)d\mu \\
             \quad +\fint_{\Sigma_t}4\pi(1-\langle \partial _r, \nu \rangle)d\mu+\int_{\Sigma_t}(1-\langle \partial _r,
             \nu
            \rangle)^2d\mu.\\
        \end{multline}
    Hence, both quantities in (i) and (ii) are uniformly bounded because
    $$\int_{\Sigma_t}O(\exp(-3r))d\mu=O(\exp(2\overline{r}_t-3\underline{r}_t)).$$
    Integrating the identity
        $$\Delta(\exp(-jr)/j)=-\exp(-jr)\Delta r+j\exp(-jr)\dr^2$$
    we obtain, after using Lemma \ref{h} and Proposition \ref{laplace},
        \begin{multline*}
            \int_{\Sigma_t}j\dr^2\exp(-jr)d\mu={\int}_{\Sigma_t}(4-2\dr^2)\exp(-(j+2)r)d\mu\\
            -4\pi{\fint}_{\Sigma_t}\exp(-jr)d\mu+{\fint}_{\Sigma_t}4\pi(1-\langle\partial_r,\nu\rangle)\exp(-jr)d\mu\\
             +\int_{\Sigma_t}(1-\langle\partial_r,\nu\rangle)^2\exp(-jr)d\mu+\int_{\Sigma_t}O(\exp(-(3+j)r))d\mu\\
             +\fint_{\Sigma_t}O(\exp(-3r))d\mu\int_{\Sigma_t}O(\exp(-3r))d\mu.
        \end{multline*}
    As a result, the third identity follows. Next, we use (iii) and the Euclidean isoperimetric identity  to prove both
    (i) and (ii).

    With respect to the the unit ball model for the hyperbolic metric, the metric $g$ can be written  as
        $$g=\psi^2\left(dx^2+dy^2+dz^2\right)+O(\psi^{-3}),$$
    where
        $$\psi(x)\equiv \frac{2}{1-\lvert x\rvert^2}\;\mbox{ and }\;\sinh r=\frac{2\lvert x\rvert}{1-\lvert x\rvert^2}.$$
    If $\underline{s}_t$ denotes the radius of the largest Euclidean ball centered at the origin that is
    contained in the interior of $\Sigma_t$, then $\sinh\underline{r}_t=2\underline{s_t}/(1-\underline{s_t}^2)$.
    The Euclidean
    isoperimetric inequality implies that
        $$\h2(\Sigma_t)\geq (36\pi)^{1/3}\v(\{\lvert x\rvert = \underline{s}_t\})^{2/3}=4\pi\underline{s}_t^2,$$
    where the area and volume are measured with respect to the Euclidean metric. Thus, denoting by $d\sigma$ the
    surface measure induced by the Euclidean metric, we have
        \begin{multline}\label{24}
            \int_{\Sigma_t}4\exp(-2r)d\mu \geq \int_{\Sigma_t}\sinh^{-2}rd\mu-\int_{\Sigma_t}\exp(-2r)\sinh^{-2} rd\mu\\
         = \int_{\Sigma_t}\lvert x \rvert^{-2} d\sigma+O(\exp(-2\underline{r}_t))
         \geq 4\pi-4\pi(1-\underline{s}^2)+O(\exp(-2\underline{r}_t))\\
         = 4\pi+O(\exp(-\underline{r}_t)).
        \end{multline}

    On the other hand,
        $$\fint_{\Sigma_t}(1-\langle \partial _r, \nu \rangle)d\mu\leq \a_t^{-1/2}\left(\int_{\Sigma_t}
        (1-\langle\partial_r,\nu\rangle)^2d\mu\right)^{1/2}=O(\exp(-2\underline{r}_t))$$
    and, recalling that both quantities in (i) and (ii) are uniformly bounded, we obtain from \eqref{e1} that
    $$
        4\pi+O(\exp(-\underline{r}_t))={\int}_{\Sigma_t}4\exp{(-2r)}d\mu+\int_{\Sigma_t}(1-\langle\partial_r,\nu\rangle)^2d\mu.
    $$
    Therefore, (i) and (ii) follow  from this estimate combined with \eqref{24}.
\end{proof}

Next, we use the stability of $\Sigma_t$ in the same way as in \cite[Section 5]{Huisken3} in order to obtain
integral estimates for $\mathring{A}$.

\begin{prop}\label{integral}
The following estimate holds
    \begin{multline*}
            \fint_{\Sigma_t}\n^2d\mu+\int_{\Sigma_t}\n^4d\mu+\int_{\Sigma_t}\lvert\nabla
            \tf\rvert^2d\mu\leq\int_{\Sigma_t}\dr^2O(\exp(-4r))d\mu\\
            +\int_{\Sigma_t}O(\exp(-7{r}))d\mu.
    \end{multline*}
In particular,
$$\fint_{\Sigma_t}\n^2d\mu\leq O(\exp(-4\underline{r}_t)). $$
\end{prop}
\begin{proof}
    Integrating the identity in Lemma \ref{A} and using Lemma \ref{h} we obtain
        \begin{multline*}
            8\pi\fint_{\Sigma_t}\n^2d\mu+\int_{\Sigma_t}\lvert\nabla \tf\rvert^2d\mu=\int_{\Sigma_t}\left(\n^2+\n\dr^2\right)O(\exp(-3r))d\mu\\
            +\int_{\Sigma_t}\n^4d\mu+\int_{\Sigma_t}\n^2d\mu\fint_{\Sigma_t}O(\exp(-3r))d\mu+\int_{\Sigma_t}\n O(\exp(-5r))
            d\mu.
        \end{multline*}
    We now argue that, for every fixed $\varepsilon>0,$
    \begin{multline}\label{f1}
        (8\pi+o(1))\fint_{\Sigma_t}\n^2d\mu+\int_{\Sigma_t}\lvert\nabla \tf\rvert^2d\mu=(1+\varepsilon/2)\int_{\Sigma_t}\n^4d\mu\\
        +\int_{\Sigma_t}\dr^2O(\exp(-4r))d\mu+\int_{\Sigma_t}O(\exp(-7{r}))d\mu.
    \end{multline}
    This is true because
    $$ \int_{\Sigma_t}\n O(\exp(-5{r}))d\mu\leq \int_{\Sigma_t}\n^2 O(\exp(-3\underline{r}_t))d\mu+
    \int_{\Sigma_t} O(\exp(-7{r}))d\mu,$$
    $$\n\dr^2 O(\exp(-3r))\leq \varepsilon/2\n^4+\dr^2 O(\exp(-4r)),$$
    and, due to Proposition \ref{estimates},
        $$\lim_{t\to\infty} \fint_{\Sigma_t}O(\exp(-3r))d\mu=0.$$

    Before we use the stability of $\Sigma_t$ let us first remark that, according to Lemma \ref{ads} and Lemma
    \ref{h}, we have for all $t$ sufficiently large
        $$\lvert A \rvert^2 +R(\nu,\nu) \geq   \n^2+2+R(\nu,\nu) \geq \n^2+O(\exp(-3r)).$$
    Hence, the stability assumption implies that
            $$\int_{\Sigma_t}\n^2 f^2d\mu\leq \int_{\Sigma_t}\lvert \nabla f \rvert^2 d\mu+\int_{\Sigma_t}f^2O(\exp(-3r))d\mu$$
    for all functions $f$ with $\int_{\Sigma_t}fd\mu=0.$
    Using the test function (see \cite[Section 5]{Huisken3})
        $$f=u-\fint_{\Sigma_t}ud\mu\equiv u-\bar{u},\quad\mbox{with}\quad u=\n,$$
    we obtain
        \begin{multline*}
            \int_{\Sigma_t}\n^4 d\mu \leq \int_{\Sigma_t}\lvert\nabla\n\rvert^2d\mu+2\bar{u}\int_{\Sigma_t}\n^3 d\mu\\
            +\int_{\Sigma_t}(u-\bar{u})^2O(\exp(-3r)) d\mu.
        \end{multline*}
    Looking at the proof of Lemma \ref{h} we see that
        $$\int_{\Sigma_t}\n^2d\mu=\int_{\Sigma_t}O(\exp(-3r))d\mu=o(1)$$
    and so, because
        $$ \bar{u}^2\leq\fint_{\Sigma_t}\n^2d\mu,$$
    we obtain for every $\varepsilon>0$ fixed
        \begin{multline}\label{f2}
            \int_{\Sigma_t}\n^4d\mu\leq \varepsilon\int_{\Sigma_t}\n^4d\mu+o(1)\fint_{\Sigma_t}\n^2d\mu\\
            +\int_{\Sigma_t}\lvert\nabla\n\rvert^2d\mu+\int_{\Sigma_t}(u-\bar{u})^2O(\exp(-3r))d\mu.
        \end{multline}
    We now estimate the last two terms in this inequality. If $\{e_1,e_2\}$ denotes an eigenbasis for $\tf$, one can
    easily check that
        $$ \lvert\nabla\n\rvert^2=2\lvert\nabla\tf(e_1,e_1)\rvert^2, \quad\nabla\tf(e_1,e_1)=-\nabla\tf(e_2,e_2)$$
        and
        $$\lvert\nabla\tf\rvert^2=\lvert\nabla\n\rvert^2+2\lvert\nabla\tf(e_1,e_2)\rvert^2.$$
    Using Codazzi equations and Lemma \ref{ads} we have that, for $i,j=1,2,$
        \begin{align*}
            \nabla_{i}\tf(e_i,e_j)&=\nabla_{j}\tf(e_i,e_i)+R(\nu,e_j)\\
        &=\nabla_{j}\tf(e_i,e_i)+\dr O(\exp(-3r))+O(\exp(-5r))
    \end{align*}
    and this implies that, for every $\varepsilon>0$,
        \begin{multline*}
            |\nabla_{i}\tf(e_i,e_j)|^2\geq(1-\varepsilon/2)|\nabla_{j}\tf(e_i,e_i)|^2\\+\dr^2 O(\exp(-6r))+O(\exp(-10r)).
        \end{multline*}
    Therefore, we can estimate
      $$2\lvert\nabla\tf(e_1,e_2)\rvert^2 \geq (1-\varepsilon)\lvert\nabla\n\rvert^2+\dr^2 O(\exp(-6r))+O(\exp(-10r))$$
    and obtain that
        \begin{multline*}
        \int_{\Sigma_t} \lvert\nabla\tf\rvert^2 d\mu \geq (2-\varepsilon)\int_{\Sigma_t}\lvert\nabla\n\rvert^2d\mu\\+\int_{\Sigma_t}\dr^2 O(\exp(-6r))d\mu+\int_{\Sigma_t}O(\exp(-10r))d\mu.
    \end{multline*}
    To handle the last term we remark that
        $$\int_{\Sigma_t}(u-\bar{u})^2O(\exp(-3r))d\mu\leq o(1)\fint_{\Sigma_t}\n^2d\mu $$
    and hence, we can rewrite equation \eqref{f2} as
        \begin{multline*}
            (1-\varepsilon)\int_{\Sigma_t}\n^4d\mu\leq
            1/(2-\varepsilon)\int_{\Sigma_t}\lvert\nabla\tf\rvert^2d\mu+o(1)\fint_{\Sigma_t}\n^2d\mu\\
            +\int_{\Sigma_t}\dr^2 O(\exp(-6r))d\mu+\int_{\Sigma_t}O(\exp(-10r))d\mu.
        \end{multline*}
    Multiplying this inequality by $(1+\varepsilon)/(1-\varepsilon)$ (with $\varepsilon$ small) and adding to equation \eqref{f1} we obtain
        \begin{multline*}
            \fint_{\Sigma_t}\n^2d\mu+\int_{\Sigma_t}\n^4d\mu+\int_{\Sigma_t}\lvert\nabla \tf\rvert^2d\mu\leq\int_{\Sigma_t}\dr^2O(\exp(-4r))d\mu\\+\int_{\Sigma_t}O(\exp(-7{r}))d\mu.
        \end{multline*}
\end{proof}

\section{Intrinsic geometry}\label{Intrin}

We study the intrinsic geometry of $(\Sigma_t)_{t\geq 0}$, the foliation of stable spheres with constant mean
curvature satisfying condition \eqref{condi}. More precisely, we show

\begin{thm}\label{intrinsic}
After pulling back by a suitable diffeomorphism from $\Sigma_t$ to $S^2$, the metric
    $$\hat{g}_t\equiv 4\pi\a_t^{-1}{g_t}$$
can be written as $$\exp(2\beta_t)g_0$$ with
    $$\sup\lvert\beta_t\rvert=O(\exp(2\overline{r}_t-3\underline{r}_t))\quad\mbox{and}\quad
    \int_{S^2}|\nabla \beta_t|^2d\mu_0=O(\exp(4\overline{r}_t-6\underline{r}_t)),$$
     where the norms are computed with respect to $g_0$, the standard round metric on $S^2$.
\end{thm}

We need to show that the Gaussian curvature $\widehat K_t$ of $\Sigma_t$ (computed with respect to $\hat g_t$)
converges to one with order $O(\exp(2\overline{r}_t-3\underline{r}_t))$. In order to do so, we know from Gauss
equation (see Lemma \ref{ricci}) that we need to estimate $\n$.
\begin{lemm}\label{beta}
 $\n^2$ is uniformly bounded.
\end{lemm}
\begin{proof}
    Suppose there is a sequence $t_i$ going to infinity and  a sequence of points $x_i$ in each $\Sigma_{t_i}$ (denoted simply by $\Sigma_i$) such that
        $$\sup_{\Sigma_i}\n=\n(x_i)\equiv1/\sigma_i\quad\mbox{ and }\quad\lim_{i\to\infty}\n(x_i)=\infty.$$
    Consider the sequence of ambient metrics ${g}_i=\sigma_i^{-2}\,g$ and denote the various geometric quantities
    with respect to $g_i$ using an index $i$. Because $H_i=\sigma_i H$ converges to zero (see Lemma \ref{h}), we have
     that for all $i$ sufficiently large
     $${\lvert A\rvert}^2_i=H_i^2/2+\n^2_i=\sigma_i^2 (H^2/2+\n^2)\leq 2.$$
    Thus, there exists a universal constant $r_0$ for which ${B}^i_{r_0}(x_i)\cap\Sigma_i$ is the
    graph over $T_{x_i}\Sigma_i$ of a function with gradient bounded by one (see, for instance, \cite{Colding})). As a result, there is a uniform
    constant $C$ such that, for all $s\leq r_0$,
        $$C^{-1}s^2\leq\int_{{B}^i_s(x_i)\cap\Sigma_i}d\mu_i\leq C s^2.$$

    Furthermore, the generalization of Michael-Simon Sobolev inequality proven in \cite{Hoffman} states the
    existence of some other universal constant $C$ such that, for every compactly supported function $u$,
        $$\left({\int_{\Sigma_i}u^2\,d\mu_i}\right)^{1/2}\leq C\left(\int_{\Sigma_i}\lvert{\nabla}u\rvert_i
        \,d\mu_i+\int_{\Sigma_i}{H}_i\lvert u\rvert \,d\mu_i\right). $$
    Hence, because $H_i$ converges to zero, we have that for all $i$  sufficiently large and every compactly supported function $u$
        $$\left({\int_{{B}^i_{r_0}(x_i)\cap\Sigma_i}u^2\,d\mu_i}\right)^{1/2}\leq
        2C\int_{{B}^i_{r_0}(x_i)\cap\Sigma_i}\lvert{\nabla}u\rvert_i\,d\mu_i.$$
    Finally, because $\n_i$ is uniformly bounded, it follows easily from Lemma \ref{A} that
        $$\Delta_i\n_i^2\geq -3\n_i^2+O(\exp(-3r)).$$
    We have now all the necessary conditions to apply Moser's iteration argument (see, for instance, \cite[Lemma 11.1.]{Li}) and obtain that, for some
    constant $C$,
        \begin{multline*}
            1=\n_i^2(x_i)\leq C \int_{\Sigma_i}\n_i^2d\mu_i+O(\exp(-3\underline{r}_t))\\
            =C\int_{\Sigma_{t_i}}\n^2d\mu+O(\exp(-3\underline{r}_t)).
        \end{multline*}
    The last expression converges to zero by Proposition \ref{integral} and this gives us a contradiction.
\end{proof}

Next, we improve the estimate on $\n^2$. The idea is to exploit the fact that the hyperbolic metric is conformal
to the Euclidean metric on the unit ball. Therefore,  each $\Sigma_t$ inherits another induced metric. The proof
consists in showing that, for this new induced metric,  we can apply Moser's iteration argument for $\n^2$ and
thus bound the supremum  by its $L^2$ norm (computed with respect to the Euclidean metric).

\begin{prop}\label{a}The following estimate holds
        $$\sup_{\Sigma_t}\n^2\leq O(\exp(2\overline{r}_t-2\underline{r}_t))\int_{\Sigma_t}\dr^2O(\exp(-4r))d\mu
        +O(\exp(2\overline{r}_t-7\underline{r}_t)).$$
        In particular,
        $$\sup_{\Sigma_t}\n^2\leq O(\exp(2\overline{r}_t-6\underline{r}_t)).$$
\end{prop}
\begin{proof}
Recall  that, with respect to the the unit model for the hyperbolic metric, the metric $g$ can be written  as
    $$g=\psi^2\left(dx^2+dy^2+dz^2\right)+O(\psi^{-3})$$
where
    $$\psi(x)=\frac{2}{1-\lvert x\rvert^2}\;\mbox{ and }\;\sinh r=\frac{2\lvert x\rvert}{1-\lvert x\rvert^2}.$$


The surfaces $\Sigma_t$ converge pointwise to the sphere of radius one and, according to Proposition,
\ref{estimates} (i)
    $$\lim_{t \to \infty}\int_{\Sigma_t}\exp{(-2r)}d\mu=\lim_{t\to\infty}\int_{\Sigma_t}d\sigma=\pi,$$
where $d\sigma$ denotes the surface measure induced by the Euclidean metric. The trace free part of the second
fundamental form and the mean curvature with respect to the Euclidean metric are denoted by $\overline{\n}$ and
$\overline{H}$ respectively. Next, we argue that we have the necessary conditions to apply Moser's iteration on
$\Sigma_t$ with respect to the Euclidean metric.

We start with the remark that
    $$ \lim_{t\to\infty}\int_{\Sigma_t}\overline{\n}^2 d\sigma=0.$$
The reason is that, by conformal invariance of $$\int_{\Sigma_t}\n^2d\mu,$$ the above identity is true if we
substitute the Euclidean metric by
 $\delta\equiv\psi^{-2}g$. Nevertheless, this metric is a lower order perturbation of the Euclidean metric and so, using the
formulas derived in \cite[Section 7]{Huisken2} (more precisely identity (7.10)) the remark follows. Therefore,
Gauss-Bonnet Theorem implies that
    $$\lim_{t\to\infty}\int_{\Sigma_t}\overline{H}^2 d\sigma=16\pi.$$
We use this to argue that the mean curvature has no concentration points.
\begin{lemm}\label{concentr}
 For every $\varepsilon_0>0$
there is a positive $r_0$ such that, for all $t$ sufficiently large and all $x$ in  $\Sigma_t$,
$$\int_{\Sigma_t\cap\overline{B}_{r_0}(x)}\overline{H}^2 d\sigma\leq \varepsilon_0.$$
\end{lemm}
\begin{proof}
Denote by $\mathcal{H}^2$ the 2-dimensional Hausdorff measure and consider the sequence of Radon measures given
by
    $$\eta_i(A)\equiv \int_{\Sigma_t\cap A}\overline{H}^2d\sigma,$$
 where $A$ is any $\mathcal{H}^2$-measurable set.  From   Allard's compactness Theorem we can extract a sequence
 $\Sigma_i$ that converges to the
unit sphere in the varifold sense and, moreover, we can also assume that $\eta_i$ converges to a Radon measure
$\bar\eta$ supported on the unit sphere. Lower semicontinuity implies that, for every $\mathcal{H}^2$-measurable
set $A$,
\begin{equation}\label{hausdorff}
    \bar\eta(A)\geq 4\mathcal{H}^2(A\cap\{\lvert x \rvert=1\}).
\end{equation}
On the other hand,
    $$4\mathcal{H}^2(\{\lvert x \rvert=1\})=16\pi=\lim_{i\to\infty} \int_{\Sigma_i}\overline{H}^2 d\sigma=
    \bar\eta(\{\lvert x \rvert=1\}).$$
Therefore, an equality in \eqref{hausdorff} must hold for every $\mathcal{H}^2$-measurable set $A$ and this
implies the desired result.
\end{proof}

An immediate consequence of this lemma combined with the generalization of  Michael-Simon Sobolev inequality
proven in \cite{Hoffman} is the existence of  universal constants $C$ and $r_0$ such that, for every Euclidean
ball $\bar{B}_{r_0}$ centered at a point in  $\{\lvert x \rvert=1\}$ and for all $t$ sufficiently large, we have
both that
    $$\left({\int_{\Sigma_t\cap \bar B_{r_0}}u^2\,d\sigma}\right)^{1/2}\leq C\int_{\Sigma_t\cap \bar B_{r_0}}\lvert\overline{\nabla}u\rvert\,d\sigma$$
for every compactly supported function $u$ and that
    $$\sqrt{\mathcal{H}^2(A)}\leq K\mathcal{H}^1(\partial A)$$
for every open subset $A$ of $\Sigma_t\cap \bar B_{r_0}$ with rectifiable boundary. A standard argument implies
the existence of a universal constant $C$ for which
    $$C^{-1}s^2\leq\mathcal{H}^2\left(\widehat{B}_s(x)\right)\leq Cs^2\quad\mbox{ for all }\quad s\leq r_0,$$
where $x$ is in $\Sigma_t$ and $\widehat{B}_s(x)$ denotes the intrinsic ball of radius $s$.

With respect to the Euclidean Laplacian, the equation for $\n^2$ becomes (see Lemma \ref{A})

\begin{multline*}
     \overline{\Delta}{\lvert\mathring{A}\rvert^2} \geq  -4\psi^{2}\lvert \tf \rvert ^2 \lvert \tf \rvert^2+O(1)\n^2
     +\dr^2\n O(\exp(-r))\\+O(\exp(-6r)).
\end{multline*}
 The zeroth order terms are bounded and so, in order to apply Moser's iteration argument, we need to check that
  $\psi^2\n^2$ is in $L^{p}$ for some $p>1$. Because $\n$ is bounded, we have from Proposition \ref{integral}
    \begin{align*}
        \int_{\Sigma_t}\psi^{2+\varepsilon}\n^{2+\varepsilon}d\sigma &
        =\int_{\Sigma_t}\psi^{\varepsilon}\n^{2+\varepsilon}d\mu+
        \int_{\Sigma_t}\n^{2+\varepsilon}O(\exp(-(3-\varepsilon)r))d\mu\\
        & = O(\exp((2+\varepsilon)\overline{r}_t-4\underline{r}_t)),
    \end{align*}
which is bounded for all $\varepsilon$ sufficiently small. Hence, setting
    $$\alpha_t\equiv\sup_{\Sigma_t}\n,$$
Moser's iteration (see, for instance, \cite[Lemma 11.1]{Li}) implies the existence of some constant $C$
depending on $\varepsilon$ small for which
\begin{multline*}
 \alpha_t^2\leq C\int_{\Sigma_t}\n^2d\sigma+\alpha_t\left(
\int_{\Sigma_t}\dr^{2+2\varepsilon}O(\exp(-(1+\varepsilon)r))
d\sigma\right)^{1/(1+\varepsilon)}\\+O(\exp(-6\underline{r}_t)).
\end{multline*}
Therefore, Proposition \ref{estimates} (iii) implies that
$$\alpha_t^2\leq O(\exp(-2\underline{r}_t))\int_{\Sigma_t}\n^2 d\mu+O(\exp(-2(3+\varepsilon)/(1+\varepsilon)\underline{r}_t)).$$
The result follows from applying Proposition \ref{integral} and choosing $\varepsilon$ appropriately small.

\end{proof}

We can now prove Theorem \ref{intrinsic}.

\begin{proof}[{\bf Proof of Theorem \ref{intrinsic}}]
    Lemma \ref{ricci} and Lemma \ref{h} imply that the Gaussian curvature of $\Sigma_t$ with respect to $\hat
    g_t$ satisfies
        \begin{equation*}
            \widehat{K}_t =1-\a_t\n^2/(8\pi)+O(\exp(2\underline{r}_t-3\underline{r}_t)).
        \end{equation*}
    As a result, we obtain from Proposition \ref{a} that
    $$\widehat{K}_t=1+O(\exp(4\overline{r}_t-6\underline{r}_t))+O(\exp(2\overline{r}_t-3\underline{r}_t))
    =1+O(\exp(2\overline{r}_t-3\underline{r}_t)).$$
    Hence, because we are assuming condition \eqref{condi}, $\widehat{K}_t$ converges uniformly to one and this implies
    that, after pulling back by a diffeomorphism, the metric
    $\hat{g}_t$ can be written as $\exp(2\beta_t)g_0$ where, according to \cite[Lemma 3.7]{AlicePaul},
    $\beta_t$ converges uniformly to zero and, denoting the coordinate function on $S^2$ by $x_1,x_2,$ and
    $x_3$,
    \begin{equation}\label{c1}
       \int_{S^2}x_j \exp{(2\beta_t)}d{\mu}_0 =0\quad\mbox{for }\,j=1,2,3.
    \end{equation}

    Using the smallness of  $\beta_t$ for $t$ sufficiently large, we have that $\beta_t$ satisfies the equation
    \begin{align*}
        \Delta_0 \beta_t & = 1-\widehat{K}_t\exp(2\beta_t)\\
                         & =1-\exp(2\beta_t)+O(\exp(2\overline{r}_t-3\underline{r}_t))\exp(2\beta_t)\\
                         & = -2 \beta_t+f(\beta_t)+O(\exp(2\overline{r}_t-3\underline{r}_t))\exp(2\beta_t),
    \end{align*}
    where $$f(s)\equiv 1+2s-\exp(2s)$$ is such that $\lvert f(\beta_t)\rvert= O(\beta_t^2)$. Therefore,
    \begin{multline}\label{curtis}
    \bar\beta_t\equiv\fint_{S^2}\beta_td\mu_0=\fint_{S^2}f(\beta_t)/2d\mu_0+O(\exp(2\overline{r}_t-3\underline{r}_t))\\=
    O(|\beta_t|^2_2)+O(\exp(2\overline{r}_t-3\underline{r}_t))
    \end{multline}
    and, combining integration by parts  with  Cauchy's inequality,
    \begin{equation}\label{c2}
        \int_{S^2}\lvert\nabla \beta_t\rvert^2d\mu_0\leq (2+o(1)+1/2)\int_{S^2}\beta^2_t
        d\mu_0+O(\exp(4\overline{r}_t-6\underline{r}_t)).
    \end{equation}
    On the other hand, we know that the $L^2$ norm of the projection of $\beta_t$ on the kernel of $\Delta_0+2$ has order
    $O(\lvert\beta_t\rvert^2_2)$ because, from \eqref{c1}, we have that for $j=1,2,3$,
        $$ 2\int_{S^2}x_j \beta_t d{\mu}_0 = \int_{S^2}x_j f(\beta_t) d{\mu}_0=O(\lvert\beta_t\rvert^2_2).$$
    Hence, because $\beta_t$ converges uniformly to zero,
    $$\int_{S^2}\lvert\nabla \beta_t\rvert^2d\mu_0\geq (6+o(1))
    \int_{S^2}(\beta-\bar\beta_t)^2 d\mu_0.$$
    Combining this with \eqref{curtis} and \eqref{c2} we obtain
        $$\int_{S^2}\beta^2_t d\mu_0=O(\exp(4\overline{r}_t-6\underline{r}_t))$$
    and so, by \eqref{c2}, the estimate on the $L^2$-norm of the gradient follows.
    Finally, a simple computation shows that for some constant $C$
    \begin{align*}
        \Delta_0 \beta_t^2 &\geq -C\beta_t^2-O(\exp(4\overline{r}_t-6\underline{r}_t)).
    \end{align*}
    Thus, from Moser's iteration, we obtain that for another constant $C$
    $$ \sup_{\Sigma_t} \beta_t^2\leq C\int_{S^2}\beta^2_t d\mu_0+O(\exp(4\overline{r}_t-6\underline{r}_t))
    =O(\exp(4\overline{r}_t-6\underline{r}_t)).$$

\end{proof}

\section{Unique approximation to coordinate spheres}\label{Unique}

The purpose of this section is to show that a foliation by stable spheres $(\Sigma_t)_{t\geq 0}$ with constant
mean
curvature must approach  the coordinate spheres $\{\lvert x \rvert=r\}$ when $t$ goes to infinity. 
Such result is obviously false in hyperbolic space because we can always apply an isometry to our foliation that
changes its center. Hence, we need to use the fact that the ambient manifold is an asymptotically hyperbolic
space with nonzero  mass. The connection between nonzero mass and uniqueness of the limit for any foliation will
be made through the Kazdan-Warner identity. More precisely, we show
\begin{thm}\label{unique}
Let  $(\Sigma_t)_{t\geq 0}$ be a foliation by stable spheres with constant mean curvature satisfying condition
\eqref{condi}. Then
    $$\lim_{t\to\infty}(\overline{r}_t-\underline{r}_t)=0\quad\mbox{ and }\quad\lim_{t\to\infty}\int_{\Sigma_t}\dr^2d\mu=0. $$
\end{thm}

\begin{proof}
Let
    $$w_t(x)\equiv r(x)-\hat{r}_t,\quad\mbox{ where }\quad \a_t=4\pi\sinh^2\hat{r}_t.$$
We want to show that
    $$\lim_{t\to\infty}w_t=0\quad\mbox{ and }\quad\lim_{t\to\infty}\int_{\Sigma_t}\lvert\nabla w_t \rvert^2d\mu=0.$$
An immediate consequence of Proposition \ref{estimates} and Theorem \ref{intrinsic} is that
\begin{multline}\label{mayfield}
    \lim_{t\to\infty}\int_{S^2}\exp{(-2w_t)}d{\mu}_0=\lim_{t\to\infty}\int_{\Sigma_t}4\pi\exp{(-2w_t)}/|\Sigma_t|d{{\mu}}
    \\=\lim_{t\to\infty}\int_{\Sigma_t}4\exp{(-2r)}d{{\mu}}=4\pi.
\end{multline}

We start by deriving a sequence of preliminary results. Combining the first identity in Proposition
\ref{laplace} with both Lemma \ref{h} and Theorem \ref{intrinsic},  we obtain that the Laplacian of $w_t$ with
respect to the standard round metric ${g}_0$ is given by
    \begin{equation}\label{e2}
        {\Delta}_0 w_t = \exp(-2w_t)- 1+P_t,
    \end{equation}
where
    \begin{multline*}
        P_t=  -2\lvert \nabla_0 w_t\rvert^2\exp(-2r) +(4\pi)^{-1}\a_t\exp(2\beta_t)(1-\langle \partial _r, \nu \rangle)^2\\
        +\exp(2\beta_t)(1-\langle \partial _r, \nu \rangle)
        + (\exp(-2w_t)+1)O(\exp(2\overline{r}_t-3\underline{r}_t))\\+O(\exp(-2\underline{r}_t)).
    \end{multline*}

Furthermore,  Proposition \ref{estimates} and \eqref{mayfield} imply
\begin{lemm}\label{pt}
$$\int_{\Sigma_t}\lvert P\rvert d{\mu}_0=O(\exp(-\underline{r}_t))+O(\exp(2\overline{r}_t-3\underline{r}_t)).$$
\end{lemm}

Note that if $\hat{g}_t$ was the round metric and   $P_t\equiv 0$, then solutions to equation \eqref{e2} would
correspond to constant scalar curvature metrics. Because there is no compactness for the set of such solutions,
we cannot expect to derive apriori estimates solely from equation \eqref{e2}. In order to have good estimates we
need to prevent  $\exp(-2w_t)$ from concentrating at any point. This will be achieved because of the lemma we
present next. It uses the Kazdan-Warner identity \cite{Kazdan} combined with the fact that the mass of the
ambient manifold is not zero and that $\Sigma_t$ has constant mean curvature.

\begin{lemm}\label{kw} For each of the standard coordinate functions $x_1,x_2,x_3$ on $S^2$ we have
    \begin{multline*}
    m\int_{S^2} x_i \exp(-3 w_t) d{\mu}_0=O(\exp(2\overline{r}_t-3\underline{r}_t))\int_{S^2} x_i \exp(-3 w_t) d{\mu}_0\\
    +O(\exp(3\overline{r}_t-4\underline{r}_t))+
    O(\exp(5\overline{r}_t-6\underline{r}_t)).
    \end{multline*}
In particular, for $i=1,2,3,$
    $$\lim_{t\to\infty}\int_{S^2}x_i\exp{(-3w_t)}d{\mu}_0\left(\int_{S^2}\exp{(-3w_t)}d{\mu}_0\right)^{-1}=0. $$
\end{lemm}
\begin{proof} Recall that the Gaussian curvature of $\Sigma_t$ with respect to $\hat g_t$ is denoted by
$\widehat{K}_t$.
 Kazdan-Warner identity \cite{Kazdan} says that, for the standard coordinate functions $x_1,x_2,x_3$ on $S^2,$
 the following identity holds for each $i=1,2,3,$
 $$\int_{S^2}\langle\nabla \widehat{K}_t,\nabla x_i\rangle\exp(2\beta_t)d{\mu_0}=0$$
 or, equivalently,
 $$\int_{S^2}x_i\widehat{K}_t\exp(2\beta_t)d{\mu_0}-\int_{S^2}\widehat K_t\langle\nabla \beta_t,\nabla x_i
 \rangle\exp(2\beta_t)d{\mu_0}=0.$$
We saw in Section \ref{Intrin} that $$\widehat K_t=1+O(\exp(2\overline{r}_t-3\underline{r}_t))$$ and this
implied that (see Theorem \ref{intrinsic})
$$ \beta_t=O(\exp(2\overline{r}_t-3\underline{r}_t))\quad\mbox{and}\quad \int_{S^2}|\nabla \beta_t|^2d\mu_0=O(\exp(4\overline{r}_t-6\underline{r}_t)).$$ Thus,
$$\int_{S^2} |\langle\nabla \beta_t,\nabla x_i
 \rangle|d{\mu_0}\leq O(\exp(2\overline{r}_t-3\underline{r}_t)) $$
and so
\begin{multline*}
\int_{S^2}\widehat K_t\langle\nabla \beta_t,\nabla x_i
 \rangle\exp(2\beta_t)d{\mu_0}=\int_{S^2}\langle\nabla \beta_t,\nabla x_i\rangle d{\mu_0}+O(\exp(4\overline{r}_t-6\underline{r}_t))\\
 =\int_{S^2}2\beta_t
x_id{\mu_0}+O(\exp(4\overline{r}_t-6\underline{r}_t))=O(\exp(4\overline{r}_t-6\underline{r}_t)),
\end{multline*}
where the last equality comes from the fact that (recall proof of Theorem \ref{intrinsic})
\begin{equation}\label{student}
       \int_{S^2}x_j \exp{(2\beta_t)}d{\mu}_0 =0\quad\mbox{for }\,j=1,2,3.
    \end{equation}
Therefore, the Kazdan-Warner identity becomes
\begin{equation}\label{brubaker}
\int_{S^2} x_i\widehat {K}_t \exp{(2\beta_t)}d{\mu_0}=O(\exp(4\overline{r}_t-6\underline{r}_t))\quad\mbox{for }i=1,2,3.
\end{equation}

From Lemma \ref{ricci} we know that
    \begin{multline*}
     4\pi\widehat{K}_t  = \a_t(H^2-4)/4+m\a_t /{\sinh^3 r}
     - 3m\dr^2\a_t/(2\sinh^3 r)\\-\n^2\a_t/2+O(\exp(2\overline{r}_t-5r))
    \end{multline*}
and hence, because $\Sigma_t$ has constant mean curvature, we obtain from \eqref{brubaker} and \eqref{student} that
\begin{multline*}
        m\int_{S^2}x_i\,\a_t^{3/2}\sinh^{-3} r\exp{(2\beta_t)} \,d{\mu_0}\\
         = \int_{\Sigma_t}3/2m\dr^2 \a_t^{3/2}\sinh^{-3} r \exp{(2\beta_t)}d\mu_0\\
        + \int_{S^2}\a_t^{3/2}\n^2 /2\exp{(2\beta_t)}\,d\mu_0\\+O(\exp(3\overline{r}_t-5\underline{r}_t))+O(\exp(5\overline{r}_t-6\underline{r}_t)).
    \end{multline*}
The first assertion in the result follows because, due to Theorem \ref{intrinsic}, Proposition \ref{estimates},
and Proposition \ref{integral},
\begin{gather*}
    \int_{\Sigma_t}3m\dr^2 \a_t^{3/2}/(2\sinh^3 r)d\mu_0=O(\exp(\overline{r}_t-3\underline{r}_t)),\\
    \int_{S^2}\n^2 \a_t^{3/2}d\mu_0=O(\exp(3\overline{r}_t-4\underline{r}_t)),
\end{gather*}
and
    \begin{multline*}
        m\int_{S^2} x_i\,\a_t^{3/2}/\sinh^3 r d{\mu}_0 =(4\pi)^{3/2}m\int_{S^2} x_i\exp(-3w_t)
        d{\mu}_0\\+O(\overline{r}_t-3\underline{r}_t)+O(-4\underline{r}_t).
    \end{multline*}
The second assertion of the lemma follows from H\"{o}lder's inequality and condition \eqref{condi}.
\end{proof}

We use this lemma to prove
\begin{prop}\label{ap} The functions $w_t$ are uniformly bounded in $W^{1,2}$.
\end{prop}
\begin{proof}
    Given a smooth function $u$ on $S^2$, denote $$S[u]\equiv\int_{S^2}\lvert\nabla_0 u\rvert^2d\mu_0-2\int_{S^2}u d\mu_0.$$
    This functional has the property that, given $F$ a conformal transformation of $S^2$, then
    $$S[u]=S[v]\quad\mbox{where}\quad F^*(\exp(-2u)g_0)=\exp(-2v)g_0.$$
    This invariance can be used in the same way as in  \cite[Proposition 4.1]{AlicePaul} in order to show
    \begin{lemm}\label{ac}
    $S[w_t]$ is  bounded independently of $t$.
    \end{lemm}
    The proof will be given in the Appendix. From integration by parts in \eqref{e2} we obtain
    \begin{equation}\label{ip1}
        \int_{S^2}\lvert\nabla_0 w_t\rvert^2d\mu_0-\int_{S^2}w_td\mu_0= -\int_{S^2}
        w_t\exp(-2w_t)d\mu_0-\int_{S^2}w_t P_td\mu_0.
    \end{equation}
    The last term on the right hand side converges to zero because Lemma \ref{pt} implies that
    \begin{align*}
        \int_{S^2}\lvert w_t P_t \rvert d\mu_0 &\leq \lvert \overline{r}_t-\underline{r}_t\rvert \left(
        O(\exp(-\underline{r}_t))+O(\exp(2\overline{r}_t-3\underline{r}_t)\right).
    \end{align*}
    Set $$V_t\equiv\int_{S^2}\exp{(-2w_t)}d{\mu}_0,$$
    which we know that it converges to $4\pi$. Using Jensen's inequality we obtain that, for all $t$
    sufficiently large,
    \begin{equation}\label{say}
        \int_{S^2}\lvert\nabla_0 w_t\rvert^2d\mu_0-\int_{S^2}w_t\,d\mu_0\leq V_t\log\left(V_t^{-1}\int_{S^2}
        \exp(-3w_t)d\mu_0\right)+1.
    \end{equation}
    We estimate the first term on the right-hand side using a slightly modified form of Aubin's inequality
    \cite[Theorem 6]{Aubin}.
    \begin{lemm}[Aubin] Given a smooth function $u$ in $S^2$ such that for all $i=1,2,3$
        $$\left\rvert\int_{S^2}x_i\exp{(-3u)}d{\mu}_0\right\rvert\left(\int_{S^2}\exp{(-3u)}d{\mu}_0\right)^{-1}\leq 1/2$$
     then, for every $\delta>0$ there is $C(\delta)$ so that
        $$4\pi\log\fint_{S^2}\exp(-3u)d\mu_0\leq \left(\frac{9}{8}+\delta\right)\int_{S^2}\lvert\nabla_0 u\rvert^2d\mu_0-3\int_{S^2}u\,d\mu_0+C(\delta).$$
    \end{lemm}
    Aubin's proof extends to this setting with obvious modifications. Lemma \ref{kw} implies that we can apply this
    result and so, combining it with \eqref{say}, we have
    $$C(\delta)\leq ((9/8+\delta)V_t/(4\pi)-1)\int_{S^2}\lvert\nabla_0 w_t\rvert^2 d\mu_0-(3V_t/(4\pi)-1)
    \int_{S^2}w_td\mu_0$$
    for some constant $C(\delta)$. Therefore, using Lemma \ref{ac}, we obtain the existence of some constant
    $C(\delta)$ for which
          $$ C(\delta)\leq ((2\delta-3/4)V_t/(4\pi)-1)\int_{S^2}w_t\,d\mu_0$$
     and thus, choosing $\delta$ sufficiently small, we get that
     $$\int_{S^2}w_t\,d\mu_0$$
     is uniformly bounded from above. We already know it is trivially bounded from below because $V_t$ converges to
     $4\pi$ and so, we can  apply Lemma \ref{ac} and conclude that
     $$\int_{S^2}\lvert\nabla_0 w_t\rvert^2d\mu_0$$ is uniformly bounded above. The result follows from Poincar\'e inequality.
     \end{proof}

    For any sequence ${t_j}$ going to infinity, we can extract a subsequence so that $w_{t_j}$ converges weakly in
    $W^{1,2}$ to $w_0$. Furthermore, we know that from Trudinger's inequality \cite{Trudinger} (see also \cite[Corollary 1.8]{Alice2}),
    $\exp{(-pw_{t_j})}$
    converges in $L^1$ to  $\exp{(-pw_{0})}$ for any $p$. Therefore, for $i=1,2,3,$
    \begin{equation}\label{jj}
        m\int_{S^2}x_i \exp(-3 w_0) d{\mu}_0=0\quad\mbox{and}\quad \int_{S^2}\exp(-2 w_0) d{\mu}_0=1.
    \end{equation}
    On the other hand, it is easy to recognize from \eqref{e2} that $w_0$ satisfies weakly
    $$\Delta_0 w_0=\exp(-2w_0)-1$$
    and so, the identities in \eqref{jj} imply that $w_0=0$, i.e., $w_t$ converges weakly in $W^{1,2}$ to zero.
    Hence, from Rellich's Theorem and integration by parts formula \eqref{ip1} we obtain
    $$\lim_{t\to\infty}\int_{S^2}\lvert\nabla_0 w_t\rvert^2d\mu_0=\lim_{t\to\infty}\int_{S^2} w_t^2d\mu_0=0.$$

    We argue next that $w_t$ converges to zero uniformly. We start with
    \begin{lemm}
        $$\lim_{t\to \infty}\sup_{\Sigma_t}\dr^2=0.$$
    \end{lemm}
    \begin{proof}
    Note that  we have from Proposition \ref{estimates}
    $$\lim_{t\to\infty}\int_{\Sigma_t}\dr^2d{\mu}=0.$$
    Arguing like in the proof of Proposition \ref{laplace} and using the fact that the norm of the second fundamental
    form of $\Sigma_i$ is bounded, it is straightforward to see that
    $$\lvert \nabla \langle \nu, \partial _r \rangle \rvert\leq C\dr+O(\exp(-3r))$$
    for some constant $C$.
    Combining Proposition \ref{laplace}  with the Bochner formula for $\lvert\nabla r\rvert^2$ we obtain that,
     for some other constant C and for all $t$ sufficiently large,
    \begin{align*}
        \Delta  \dr^2 & \geq -C\dr^2+2\lvert\nabla\nabla r\rvert^2-\lvert\nabla\nabla r\rvert \dr^2O(\exp(-2r))\\
        & \quad+\dr O(\exp(-3r))\\
        & \geq -C\dr^2+O(\exp(-6r)).
    \end{align*}
    Because the second fundamental form of $\Sigma_t$ is bounded, we can invoke the same reasons as in Lemma \ref{beta}
    in order find $r_0$ so that, for all $t$ sufficiently large and all
    $x$ in $\Sigma_t$,  we can apply Moser's iteration argument in ${B}_{r_0}(x)\cap\Sigma_t$ and conclude that
           $$\sup_{\Sigma_t}\dr^2\leq C\int_{\Sigma_t}\dr^2d\mu+O(\exp(-6r)),$$
    where $C$ is a constant independent of $t$.
    \end{proof}
    An immediate consequence of this lemma, Proposition \ref{estimates}, and the identity
    $\dr^2=1-\langle\nu,\partial_r\rangle^2+O(\exp(-5r))$ is

    \begin{lemm}\label{bound}
    For all $t$ sufficiently large $$O(\exp(-5r))+\dr^2/3\leq 1-\langle\nu,\partial_r\rangle\leq \dr^2/2+O(\exp(-5r)).$$
    \end{lemm}

    A simple computation using identity \eqref{e2} shows that
    \begin{equation}\label{say2}
        \Delta_0 w_t^2  \geq 8\pi w_t(\exp(-2w_t)-1)+2P_tw_t
    \end{equation}
    and we want to write this inequality as
    $$ \Delta_0 w_t^2\geq -Cw_t^2-\mbox{(finite sum of positive terms }f_i),$$
    where $C$ is a $t$-independent constant and each term $f_i$ has its $L^p$ norm converging to zero for some $p>1$.
    We can then apply Moser's iteration argument and conclude that, for some constant $C$,
    $$\sup w_t^2\leq C\int_{S^2}w_t^2d\mu_0+C\sum_i|f_i|_{L^p}.$$
    This implies uniform convergence of $w_t$ to zero.

    In what follows, C is a generic constant independent of $t$. The first term in \eqref{say2} can be easily estimated as
        $$ w_t(\exp(-2w_t)-1) \geq -w_t^2/2-(\exp(-2w_t)-1)^2/2,$$
    where the term with the exponential converges to zero in any $L^p$-norm due to Trudinger's inequality
    \cite{Trudinger}. Looking at the expression of $P_t\,w_t$ (see \eqref{e2}), we have that the
    the first term and the third term cause no problem because
    \begin{align*}
        -\lvert \nabla_0 w_t\rvert^2\exp(-2r)w_t& =O(1)\dr^2\exp(-2w_t)w_t+w_tO(\exp(-2\underline{r}_t))\\
        &\geq -Cw_t^2-C\dr^4\exp(-4w_t)+O(\exp(-4\underline{r}_t))
     \end{align*}
     and also
      $$(1-\langle \partial _r, \nu \rangle)w_t\geq-w_t^2/2- (1-\langle \partial _r, \nu \rangle)^2/2\geq-w_t^2/2- \dr^4/8.$$
     The same sort of estimate works for the last two terms and to handle the second term we note that, by Proposition
     \ref{estimates},
     $$
        \int_{S^2} \ge^2\dr^2d\mu_0\leq 9 \int_{\Sigma_t}(1-\langle \partial _r, \nu
        \rangle+O(\exp(-5r)))^2d\mu=O(\exp(-\underline{r}_t)).
     $$
     Hence, for some constant C,
     \begin{align*}
        \int_{S^2}(\a_t(1-\langle \partial _r, \nu \rangle)^2\lvert w_t\rvert)^{1+\varepsilon}d\mu_0 & \leq C\int_{S^2}
        \ge^{2+2\varepsilon}\dr^{2+2\varepsilon}\lvert w_t\rvert^{1+\varepsilon}d\mu_0\\
      &\leq C\int_{S^2} \ge^{2+2\varepsilon}\dr^{2-2\varepsilon}\lvert w_t\rvert^{1+\varepsilon}d\mu_0\\
      &= C\a_t^{\varepsilon}\int_{S^2} \ge^2\dr^2\lvert w_t\rvert^{1+\varepsilon}d\mu_0\\
      & \leq \lvert\overline{r}_t-\underline{r}_t\rvert^{1+\varepsilon} O(\exp(2\varepsilon\overline{r}_t-\underline{r}_t))
     \end{align*}
     and, because we are assuming condition \eqref{condi}, we obtain that, for some sufficiently small $\varepsilon$,
      $$\lim_{t\to\infty}\int_{S^2}(\a_t(1-\langle \partial _r, \nu \rangle)^2\lvert w_t\rvert)^{1+\varepsilon}d\mu_0=0.$$
\end{proof}

\section{Strong approximation to  coordinate spheres}\label{strong}

We will start by arguing that we can choose $\tilde{r}_t$ such that
 $$ H=2\cosh \tilde{r}_t/\sinh{\tilde{r}_t}-m /\sinh^3{\tilde{r}_t}+o(\exp(-4\tilde{r_t})). $$
We will then show that, with respect to $\hat{g}_t$,
 $\Sigma_t$ is $o(\exp(-\underline{r}_t))$ ``close" to $\{\lvert x\rvert=\tilde{r}_t\}$. This is exactly the rate
 we need in order to prove uniqueness for the constant mean curvature foliation. Being more precise, rename the
 family of functions $w_t$ on $\Sigma_t$  by $$w_t(x)\equiv r(x)-\tilde{r}_t.$$ This section is devoted to show
 \begin{thm}\label{close}
    With respect to the metric $\hat{g}_t$, $$\lvert w_t \rvert_{C^{2,\alpha}}=o(\exp(-\underline{r}_t)).$$
 \end{thm}
Before proving Theorem \ref{close} we need to argue that $\tilde{r}_t$ is well defined. In order to do so, the
following proposition is important
\begin{prop}\label{r}\label{decay}
The following identity holds for each $\Sigma_t:$ $$\sup_{\Sigma_t}\,\dr^2=o(\exp(-2\underline{r}_t)).$$
\end{prop}
\begin{proof}

    Direct computation implies the following gradient estimate
    \begin{lemm}\label{grad}
        $$\lvert\nabla \dr^2\rvert\leq (4\exp(-2r)-4\pi/\a_t)\dr+3\dr^3+\dr o(\exp({-2\underline{r}_t})).$$
    \end{lemm}
    \begin{proof}
        Because $g$ is a perturbation of the hyperbolic metric,  if we decompose a unit tangent vector $V$ as
            $\overline{V}+\beta\,\partial_r$, where $\overline{V}$ has no $\partial_r$ component,
        then $\lvert\beta\rvert\leq \dr+O(\exp(-3\underline{r}_t))$ and so
        \begin{equation}\label{for1}
            \lvert\langle V-\overline{V}, \partial_r^{\top} \rangle\rvert\leq \dr^3+\dr^2O(\exp(-3\underline{r}_t)).
        \end{equation}
        Denoting the connection with respect to the hyperbolic metric by $\overline{\nabla}$, we can
        estimate
        \begin{align*}
            \langle V, \nabla \dr^2 \rangle &= 2 \langle\overline{\nabla}_V \partial_r,\partial_r^{\top}\rangle-2\langle \partial_r,\nu\rangle\langle {\nabla}_V \nu,\partial_r^{\top}\rangle+\dr
            O(\exp(-3\underline{r}_t)).
        \end{align*}
    Due to Theorem \ref{unique} and Proposition \ref{a}, we know that
        $$\sup_{\Sigma_t}\n=o(\exp(-2\underline{r}_t)).$$
    Hence,
        \begin{align*}
            \langle V, \nabla \dr^2 \rangle& =2(\cosh r/\sinh r)\,\langle\overline{V} ,\partial_r^{\top}\rangle-
            \langle \partial_r,\nu\rangle \langle{V} ,\partial_r^{\top}\rangle H\\
            &\quad+o(\exp(-2\underline{r}_t))\dr\\
            & \leq (2\cosh r/\sinh r-H)\dr+H(1-\langle \partial_r,\nu\rangle)\langle\overline{V} ,\partial_r^{\top}\rangle\\
             &\quad+H\langle \partial_r,\nu\rangle\langle\overline{V}-V ,\partial_r^{\top}\rangle+o(\exp(-2\underline{r}_t))\dr\\
             &\leq (4\exp(-2r)-4\pi/\a_t)\dr+3\dr^3+\dr o(\exp({-2\underline{r}_t})),
        \end{align*}
    where in the last inequality we used Lemma \ref{h}, Lemma \ref{bound}, and \eqref{for1}.
    \end{proof}
     Choose $p$ in $\Sigma_t$ so that $\dr(p)=0$. Given $q$ in $\Sigma_t$, denote by $\gamma(s)$ the unit speed
     geodesic connecting $p$ to $q$ and let $f(s)=\dr(\gamma(s))$. From Lemma \ref{grad} we have
    \begin{align*}
        \left\lvert  \frac{d}{ds}f^2\right\rvert&\leq (\lvert4\pi/\a_t-4\exp(-2r)\rvert+ o(\exp({-2\underline{r}_t})))f+3f^3\\
        &\leq \alpha_t^2f+3f^3,
    \end{align*}
    where $$\alpha_t^2=\sup_{\Sigma_t}\lvert4\pi/\a_t-4\exp(-2r)\rvert+ o(\exp({-2\underline{r}_t})).$$ If we set
    $$V(s)\equiv\alpha_t\tan((3/4)^{1/2}\alpha_t\, s)/\sqrt{3},$$ then $$(V^2)'=\alpha_t^2V+3V^2$$
    and so, while $V$ is well defined, we have $f \leq V$ because $f(p)=0$. Moreover, $\mathrm{diam}(\hat{g}_t)$ converges
    to $2$ and therefore
        \begin{multline*}
            \lim_{t\to\infty}\alpha_t^2\mathrm{diam}^2({g}_t)= \lim_{t\to\infty}\alpha_t^2\a_t\pi^{-1}\\
            \leq \lim_{t\to\infty}\sup_{\Sigma_t}\,\lvert 4-4\a_t\pi^{-1}\exp{(-2\underline{r}_t)} \rvert=0,
        \end{multline*}
    where the last equality follows from Theorem \ref{unique}. Hence, we obtain for all $t$ sufficiently large
        $$\dr(q)=f(q)\leq V(\mathrm{diam}({g}_t)/2)\leq \alpha_t = o(\exp(-\underline{r}_t)).$$
\end{proof}

The decay for $\dr$ proved above combined with Lemma \ref{bound} allow us to write the second identity in
Proposition \ref{laplace} as
    $$\Delta r= 2\cosh r/\sinh r- m/\sinh^3 r -H +P,$$
    where $\lvert P\rvert=o(\exp(-4r))$. Thus, we obtain after integration
    $$H=\fint_{\Sigma_t}2\cosh r/\sinh r- m/\sinh^3 r d\mu+\fint o(\exp(-4r))d\mu$$
and so,  the existence of $\tilde{r}_t$ satisfying $\underline{r}_t\leq \tilde{r}_t\leq \overline{r}_t$ follows
easily.

We also remark that a careful inspection of the term $P$ shows that
    \begin{equation}\label{P}
        \lvert\nabla P\rvert=o(\exp(-5r))
    \end{equation}
because, as can be seen from Lemma \ref{bound}, Lemma \ref{grad}, and Proposition \ref{decay}, we have
    $$\lvert\nabla\langle\partial_r, \nu\rangle \rvert=o(\exp(-3r))\quad\mbox{and}\quad\lvert\nabla \dr^2\rvert
    =o(\exp(-3r)). $$

\begin{proof}[{\bf Proof of Theorem \ref{close}}]
We know from Theorem \ref{unique} that $w_t$ converges to zero uniformly and thus, the expansion for $H$ implies
    \begin{align*}
        \Delta r 
                 &= 4\exp(-2\tilde{r}_t)(\exp(-2w_t)-1)-8m\exp(-3\tilde{r}_t)(\exp(-3w_t)-1)\\ &\quad+ o(\exp(-4r))\\
                 &=-8 w_t   \exp(-2\tilde{r}_t)+f(w_t) O\left(\exp(-2r)\right)+w_t O(\exp(-3r))\\
                 & \quad  + o(\exp(-4r)),
        \end{align*}
    where
        $$f(x)\equiv\exp(-2x)-1+2x.$$
    Hence, from Theorem \ref{intrinsic}, the above identity translates to
    \begin{equation*}
        \Delta_0 w_t + (2\exp(2\hat{r}_t-2\tilde{r}_t+2\beta_t)+o(1))w_t=f(w_t)O\left(1\right)+o(\exp(-2r)).
    \end{equation*}

When $t$ goes to infinity the above operator converges to  $\Delta_0 +2$. Therefore, in order to get a good
control of $w_t$ we need to control the projection of $w_t$ on the kernel of $\Delta_0 +2$. This is achieved
using the following improvement of Lemma \ref{kw}
\begin{lemm}
    \begin{equation*}
            \int_{S^2} w_t x_id\mu_0=o(\exp(-\underline{r}_t))+\lvert w_t\rvert_2^2 O(1), \quad i=1,2,3.
    \end{equation*}
\end{lemm}
\begin{proof}
    We saw in Lemma \ref{grad} that $\n^2$ has order $o(\exp(-4\underline{r}_t))$ and so, we obtain  from Lemma
    \ref{ricci} and Proposition \ref{decay} that
        $$\widehat{K}_t=1+o(\exp(-\underline{r}_t)).$$
    This implies that (check proof of Theorem \ref{intrinsic})
        $$\exp(2\beta_t)=1+o(\exp(-\underline{r}_t))$$
    and thus, an inspection of the proof of Lemma \ref{kw} shows that
        $$m\int_{S^2}\exp(-3w_t)x_id\mu_0=o(\exp(-\underline{r}_t))\quad i=1,2,3.$$
    This implies the desired result.
\end{proof}

As a result, if we decompose $w_t$ as $u_t+v_t,$ where $u_t$ is in the kernel of $\Delta_0 +2$ and $v_t$ is
perpendicular to $u_t$, we obtain from  the previous lemma that
    \begin{equation*}
        \lvert u_t\rvert_{C^{2,\alpha}}=o(\exp(-\underline{r}_t))+ \lvert w_t\rvert_2^2O(1)
    \end{equation*}
and thus, because $u_t$ converges to zero uniformly,
    \begin{equation}\label{e6}
        \lvert u_t\rvert_{C^{2,\alpha}}= o(\exp(-\underline{r}_t))+ \lvert v_t\rvert_2^2O(1).
    \end{equation}

On the other hand, we have that $v_t$ satisfies
    \begin{equation}\label{e4}
        \Delta_0 v_t +(2+o(1))v_t=u_t o(1)+f(w_t)O\left(1\right)+o(\exp(-2r)).
    \end{equation}
Looking at the way the terms in this equation were derived and using \eqref{P} we see that
$$\lvert u_t o(1) \rvert_{C^{0,\alpha}}=o(1)\lvert u_t \rvert_{C^{0,\alpha}},\quad
\left\lvert f(w_t)O\left(1\right) \right\rvert_{C^{0,\alpha}}\leq O(1)\lvert w_t \rvert_{C^{0,\alpha}}^2,$$ and
$$\lvert o(\exp(-2r))\rvert_{C^{0,\alpha}}= o(\exp(-2\underline{r}_t)).$$
Hence, from Schauder estimates, we have that for some constant $C$
    $$
        \lvert v_t \rvert_{C^{2,\alpha}} \leq o(1)\lvert u_t \rvert_{C^{0,\alpha}}+C \lvert w_t \rvert_{C^{0,\alpha}}^2
        +C \sup\lvert v_t\rvert+o(\exp(-2\underline{r}_t))
    $$
Moreover,  $\lvert w_t\rvert_{C^{0,\alpha}}$ converges to zero (see Proposition \ref{decay}) and hence we have
from  \eqref{e6}
 $$\lvert v_t \rvert_{C^{2,\alpha}}\leq  \sup \lvert v_t\rvert+o(\exp(-r)).$$

Using the same arguments as in the proof of Theorem \ref{intrinsic}, we can see that the orthogonality condition
of $v_t$, identity \eqref{e6}, and equation \eqref{e4}, imply the estimate
    $$\int_{S^2}w^2_td\mu_0=o(\exp(-2\underline{r}_t)).$$
A simple computation shows that, for some constant $C$,
    $$\Delta_0 v_t^2\geq-Cv_t^2 -C u_t^2+o(\exp(-2\underline{r}_t))$$
and thus, using Moser's iteration we obtain
    $$\sup v_t^2\leq C\int_{S^2}w^2_td\mu_0+C\sup u^2_t+o(\exp(-2\underline{r}_t))=o(\exp(-2\underline{r}_t)).$$
\end{proof}

\section{Existence and uniqueness of constant mean curvature foliations}

\subsection{Uniqueness}\label{uni}
We are now ready to prove the main theorem.
\begin{thm} If the mass of the asymptotically Anti de Sitter- Schwarzschild  metric is nonzero, then
    any two smooth foliations by stable spheres with constant mean curvature
    $\left(\Sigma^1_t\right)_{t\geq 0}$ and $\left(\Sigma^2_t\right)_{t\geq 0}$ for which
    \begin{equation}
            \lim_{t\to\infty}(\overline{r}_t-4/3\underline{r}_t)=-\infty
    \end{equation}
     will coincide for $t$ sufficiently large.
\end{thm}
\begin{proof}
    We can reparametrize the foliations so that $H(\Sigma^1_t)=H(\Sigma^2_t)=H_t$ for all $t$ sufficiently large.
    The results in the previous section imply the existence of $\tilde{r}_t$ so that
        $$H_t=2\cosh \tilde{r}_t/\sinh{\tilde{r}_t}-m /\sinh^3{\tilde{r}_t}+o(\exp(-4\tilde{r_t}))$$
    and for $i=1,2,$ the functions
        $$w^i_t\equiv r(x)-\tilde{r}_t $$
    have $C^{2,\alpha}$ norm of order $o(\exp(-\underline{r}_t))$, where the norm is computed with respect to the
    rescaled metric $\hat{g}^i_t$.
    From Theorem \ref{intrinsic} we know that, after pulling  back $\hat{g}^1_t$ by a suitable diffeomorphism,
    the metric can be written as $\exp(2\beta_t)g_0$, where $g_0$ denotes the standard round metric on $S^2$. The following lemma improves the estimate on $\beta_t$.
    \begin{lemm}\label{scalar}
        The functions $\beta_t$ can be chosen so that $$\lvert\beta_t\rvert_{C^{2,\alpha}}=o(\exp(-\underline{r}_t)).$$
    \end{lemm}
    \begin{proof}
        A direct computation using the fact that each $\Sigma^1_t$ is the graph over $\{\lvert x
        \rvert=\tilde{r}_t\}$ of a function $w^1_t$ with $\lvert w^1_t\rvert_{C^{2,\alpha}}=o(\exp(-\underline{r}_t))$,
        reveals that
        $\n^2_{C^{0,\alpha}}=o(\exp(-2\underline{r}_t))$, where the norm is computed with respect to the metric $g_0$.
        Combining Lemma \ref{ricci} with the asymptotic expansion of $H_t$ in terms
        of $\tilde{r}_t$, we obtain that the Gaussian curvature of $\hat{g}^1_t$ can be
        written as
            \begin{align*}
                \widehat{K}_t= 1+ S_t,\quad\mbox{where}\quad\lvert S_t\rvert_{C^{0,\alpha}}=o(\exp(-\underline{r}_t)).
            \end{align*}
        The functions $\beta_t$ were chosen so that,  for each coordinate function $x_1,x_2,$ and $x_3,$
            $$ \int_{S^2}\exp(2\beta_t)x_id\mu_0=0, \quad i=1,2,3$$
        and they satisfy the equation
            \begin{align*}
                \Delta_0\beta_t & =1-\widehat{K}_t\exp(2\beta_t)\\
                & =1-\exp(2\beta_t)+\exp(2\beta_t)o(\exp(-r)).
            \end{align*}
        One can then use the smallness of $\beta_t$ and argue in the same way as in the proof of either Theorem \ref{intrinsic}
        or Theorem  \ref{close} in order to prove the lemma.
    \end{proof}

    With respect to the coordinates $(r,\omega)$ in $\R\times S^2$ we consider, for each fixed $t$, the interpolation
    surfaces
    $$
        \Sigma_{t,s}=\left\{\tilde{r}_t+(1-s)w^1_t(\omega)+s w^2_t(\omega)\,|\, \omega\mbox{ in }S^2 \right\},
        \quad1\leq s \leq 0
    $$
    and set $$F_t(s)\equiv H(\Sigma_{t,s}).$$

In what follows, we fix $t$ large and suppress the index $t$ in the notation for the sake of simplicity. From
Taylor's formula, we have
$$F(1)=F(0)+F'(0)+P,\quad P\equiv\int_0^1\int_0^1sF''(su)du ds,$$
where $$F'(0)=\Delta_{\Sigma_1} \phi+(\lvert A \rvert^2+R(\nu,\nu))_{\Sigma_1} \phi, \quad \phi\equiv w^2-w^1$$
and, for some universal constant $C$,
$$\lvert F''(s)\rvert\leq C(\lvert \nabla^2 \phi\rvert \lvert \phi \rvert+\lvert \nabla \phi \rvert^2+\lvert \phi
\rvert^2\exp(-2r))_{|\Sigma_{t,s}}.$$
 We know that $F(1)=F(0)$ and that
    \begin{align*}
        (\lvert A \rvert^2+R(\nu,\nu))_{1} & = (H^2-4)/2-m/\sinh^3 \tilde{r}+O(\exp(-4r))\\
        & = 8\pi/\lvert\Sigma^1\rvert-3m/\sinh^3 \tilde{r}+O(\exp(-4r)).    \end{align*}
Therefore, using Theorem \ref{close}, we have that, with respect to the metric $g_0$,  $\phi$ satisfies the
following equation
    \begin{multline*}
        \Delta_0 \phi+\exp(2\beta)(2-3m/\sinh{\tilde{r}}+O(\exp(-2r)))\phi\\
        +\lvert\Sigma^1\rvert(4\pi)^{-1}\exp(2\beta)P=0
    \end{multline*}
or, alternatively,
    $$ \Delta_0 \phi+(2-3 m/\sinh{\tilde{r}})\phi=Q-\lvert\Sigma^1\rvert(4\pi)^{-1}\exp(2\beta)P,$$
where
    \begin{equation*}
        Q\equiv(1-\exp(2\beta))(2-3 m/\sinh{\tilde{r}})\phi -\exp(2\beta)O(\exp(-2r))\phi.
    \end{equation*}
Consider the decomposition $$\phi=\phi_{0}+\phi_{1},\quad Q=Q_0+Q_1,\quad
\lvert\Sigma^1\rvert\exp(2\beta)P=P_0+P_1$$ where $\phi_{0}, Q_0, P_0$ belong to the kernel of $\Delta_0+2$ and
$\phi_{1}, Q_1, P_1$ are orthogonal to $\phi_{0}, Q_0,$ $P_0$ respectively. Then
    \begin{equation}\label{e7}
         \Delta_0 \phi_1+(2-3 m/\sinh{\tilde{r}})\phi_1  =Q_1-P_1,
    \end{equation}
    \begin{equation}\label{e8}
         3m/\sinh{\tilde{r}}\,\phi_0=P_0-Q_0,
    \end{equation}
and it is immediate to recognize that, for some universal constant $C$,
$$\lvert P_0\rvert_{C^{0,\alpha}}+ \lvert P_1\rvert_{C^{0,\alpha}}\leq C \lvert  \phi \rvert^2_{C^{2,\alpha}}$$
and
$$\lvert Q_0\rvert_{C^{0,\alpha}}+ \lvert Q_1\rvert_{C^{0,\alpha}}\leq o(\exp(-\underline{r}))\lvert \phi
\rvert_{C^{0,\alpha}},$$
 where the norms are computed with respect to $g_0$.

 Applying Schauder estimates to equation \eqref{e7} and
 using the fact that $|\phi|_{C^{2,\alpha}}$ converges to zero, we obtain
$$\lvert  \phi_1 \rvert_{C^{2,\alpha}} \leq C\sup_{\Sigma^1}|\phi_1|+C\lvert  \phi _0\rvert^2_{C^{2,\alpha}}+
o(\exp(-\underline{r}))\lvert \phi_0 \rvert_{C^{0,\alpha}},$$ where $C$ is some uniform constant. The
orthogonality condition satisfied by $\phi_1$ can be used in the same way as in the proof of Theorem \ref{close}
in order to show that
$$\sup_{\Sigma^1}|\phi_1|\leq C|\phi|^2_{C^{2,\alpha}}+o(\exp(-\underline{r}))|\phi|_{C^{0,\alpha}}.$$
Therefore,
$$\lvert  \phi_1 \rvert_{C^{2,\alpha}} \leq C\lvert  \phi _0\rvert^2_{C^{2,\alpha}}+
o(\exp(-\underline{r}))\lvert \phi_0 \rvert_{C^{0,\alpha}},$$
 where $C$ is some uniform constant. Because  the ${C^{2,\alpha}}$ norm of $\phi_0$ is bounded by its
 ${C^{0,\alpha}}$ norm, we obtain from \eqref{e8} that
\begin{align*}
    \lvert  \phi_0 \rvert_{C^{2,\alpha}}& \leq O(\exp({\underline r}))\lvert  \phi \rvert^2_{C^{2,\alpha}}+o(1)\lvert
    \phi \rvert_{C^{0,\alpha}}\\
&  \leq O(\exp({\underline r}))\lvert  \phi_0 \rvert^2_{C^{2,\alpha}}+o(1)\lvert  \phi_0 \rvert_{C^{0,\alpha}}.
\end{align*}
From Theorem \ref{close} we have that $\lvert  \phi \rvert_{C^{2,\alpha}}=o(\exp(-\underline r))$  and thus
$$\lvert  \phi_0 \rvert_{C^{0,\alpha}}\leq o(1)\lvert  \phi_0 \rvert_{C^{0,\alpha}}.$$
Consequently, for $t$ sufficiently large, $\phi_0=0$ and this implies that $\phi_1=0$. Hence $w^1_t=w^2_t$ for
all $t$ sufficiently large and this is the same as $\Sigma^1_t=\Sigma^2_t$.
\end{proof}

\subsection{Existence}\label{exi}
We show existence of a foliation by stable spheres with constant mean curvature when the mass of the
asymptotically Anti de Sitter- Schwarzschild  metric is positive. This result was previously shown by Rigger \cite{Rigger} using 
a modified mean curvature flow approach. The argument we use is a straightforward adaptation of
the arguments used by Rugang Ye in \cite{Ye}, where he showed a similar theorem in the context of asymptotically
flat manifolds.  We include the proof of existence of a foliation for the sake of completeness.
\begin{thm}
If the mass of the asymptotically Anti de Sitter- Schwarzschild  metric is nonzero, the manifold $M$ admits,
outside a compact set, a foliation by spheres with constant mean curvature. Moreover, if the mass is positive,
then the spheres are stable.
\end{thm}
\begin{proof}
  Given $\phi\in C^{\infty}(S^2)$ and $r$ sufficiently large, let
  $$\Sigma_{r}(\phi)\equiv\{(r+\phi(x),x)\,|\, x\in S^2\}$$
and set $F(r,\phi)\equiv H(\Sigma_r(\phi)).$ Because the metric $g$ is a $C^3$ perturbation of $g_m$, we have
from Lemma \ref{ads} that $$F(r,0)=2\cosh r/\sinh r-m/\sinh^3 r+O(\exp(-5r))$$ and
$$(|A|^2+R(\nu,\nu))_{\Sigma_r(0)}=8\pi/|\Sigma_r(0)|-3m/\sinh^3 r+O(\exp(-5r)).$$

For every $r$ sufficiently large, we want to find $\phi$ so that
$$F(r,\phi)=2\cosh r/\sinh r-m/\sinh^3 r.$$
Using Taylor's formula like in the previous subsection, the above equation is equivalent to solve on $S^2$
$$\Delta_0 \phi+(2-3m/\sinh r)\phi=P(\phi)+Q(\phi)+N,$$
where
$$N\equiv |\Sigma_r(0)|(F(r,0)-2\cosh r/\sinh r+m/\sinh^3 r)/(4\pi)$$
satisfies $$|N|_{C^{0,\alpha}}=O(\exp(-3r))$$ and $P$ and $Q$ are such that, for some uniform constant $C$,
$$|P(\phi)|_{C^{0,\alpha}}\leq C|\phi|^2_{C^{2,\alpha}},$$
$$|P(\phi)-P(\psi)|_{C^{0,\alpha}}\leq C(|\phi|_{C^{2,\alpha}}+|\psi|_{C^{2,\alpha}})|\phi-\psi|_{C^{2,\alpha}},$$
$$|Q(\phi)|_{C^{0,\alpha}}\leq o(\exp(-r))|\phi|_{C^{0,\alpha}},$$
and
$$|Q(\phi)-Q(\psi)|_{C^{0,\alpha}}\leq o(\exp(-r))|\phi-\psi|_{C^{0,\alpha}}.$$

Consider the map $$T:C^{2,\alpha}(S^2)\longrightarrow C^{2,\alpha}(S^2)$$ such that
$$\Delta_0 (T(u))+(2-3m/\sinh r)T(u)=P(u)+Q(u)+N.$$
The map is well defined because the operator on the left hand side is invertible. The existence of a constant
mean curvature foliation follows from
\begin{lemm}\label{lee}
    For all $r$ sufficiently large, the map $T$ is a contraction of
    $$\{u\in C^{2,\alpha}(S^2)\,|\,|u|_{C^{2,\alpha}}\leq \exp(-r)/r \}$$
    onto itself.
\end{lemm}
\begin{proof}
    Set $\phi\equiv T(u)$ and consider the decomposition
        $$\phi=\phi_0+\phi_1,\quad P(u)=P_0+P_1,\quad Q(u)=Q_0+Q_1, \quad\mbox{and } N=N_0+N_1,$$
    such that $\phi_0, P_0, Q_0, N_0$ belong to the kernel of $\Delta+2$ and $\phi_1, P_1, Q_1, N_1$ are
    orthogonal to $\phi_0, P_0, Q_0, N_0$ respectively.
    Thus, $$-3m\phi_0/\sinh r=P_0+Q_0+N_0$$
    and so
    \begin{multline*}
        |\phi_0|_{C^{2,\alpha}}\leq O(\exp(r))|u|^2_{C^{2,\alpha}}+o(1)|u|_{C^{2,\alpha}}+O(\exp(-2r))\\
        \leq o(1)\exp(-r)/r.
    \end{multline*}
    Furthermore,
        $$\Delta_0 \phi_1+(2-3m/\sinh r)\phi_1=P_1(u)+Q_1(u)+N_1$$
     and hence, we can argue in the same way as in the proof of Theorem \ref{intrinsic} and
    use the orthogonality condition satisfied by $\phi_1$ in order to show that
    \begin{multline*}
        |\phi_1|_{C^{2,\alpha}}\leq C|u|^2_{C^{2,\alpha}}+o(\exp(-r))|u|_{C^{2,\alpha}}+O(\exp(-3r))\\
        \leq o(1)\exp(-r)/r
    \end{multline*}
    for some uniform constant $C$. Therefore, we have for all $r$ sufficiently large that
    $|\phi|_{C^{2,\alpha}}\leq \exp(-r)/r.$ Finally, we can argue in the same way and check that
    $$|T(u)-T(v)|_{C^{2,\alpha}}\leq o(1)|u-v|_{C^{2,\alpha}}. $$
\end{proof}

Denote by $\Sigma_r$ the constant mean curvature sphere that is the graph over $\{|x|=r\}$ of a function $w_r$
with $|w_r|_{C^{2,\alpha}}\leq \exp(-r)/r$, where the norm is computed with respect to the standard round metric
$g_0$. We need to show that $\Sigma_r$ is stable for all $r$ sufficiently large, i.e., we need to show that
second variation operator
$$ Lf \equiv-\Delta f -\left(\lvert A\rvert^2+R(\nu,\nu) \right)f$$
has only nonnegative eigenvalues when restricted to functions with zero mean value.

Because $\Sigma_r$ is the graph of a function with $|w_r|_{C^{2,\alpha}}\leq o(\exp(-r))$, we have that
$$\sup_{\Sigma_r}|A|^2=o(\exp(-4r)),$$
where this norm is computed with respect to the metric $g_r$ induced by the ambient metric. Moreover, we also
have that
$$ \lvert A\rvert^2+R(\nu,\nu)= 8\pi/|\Sigma_r|-3m/\sinh^3 r+O(\exp(-4r)).$$
In this setting, Theorem \ref{intrinsic} and Lemma \ref{scalar} apply and so, after applying a suitable
diffeomorphism, the normalized metric $\hat g_r\equiv 4\pi/|\Sigma_r|g_r$ can be written as $\exp(2\beta_r)g_0$,
where $|\beta_r|_{C^{2,\alpha}}\leq o(\exp(-r))$. Hence, in terms of the metric $g_0$, the positivity of the
operator $L$ is equivalent to the positivity of the following operator defined on $S^2$
$$L_0 \phi\equiv -\Delta_0\phi-(2-3m/\sinh r+o(\exp(-r))\phi.$$
This is true for all $r$ sufficiently large whenever $m$ is positive.
\end{proof}

\appendix
\section{Proof of Lemma \ref{ac}}
Given a smooth function $u$ on $S^2$, recall that $$S[u]\equiv\int_{S^2}\lvert\nabla_0
u\rvert^2d\mu_0-2\int_{S^2}ud\mu_0.$$ We want to show that $S[w_t]$ is bounded independently of $t$. We
essentially follow, with some necessary modifications, the proof of \cite[Proposition 4.1]{AlicePaul}. The main
idea consists in exploiting the invariance of $S$ under conformal transformations. More precisely, we will find
conformal diffeomorphisms $F_t$ for which the family of functions $u_t$ defined by
$$\exp(-2u_t)g_0=F_t^*(\exp(-2w_t)g_0)$$
is such that
$$\int_{S^2}|\nabla u_t|^2d\mu_0\quad\mbox{and}\quad\int_{S^2} u_t d\mu_0$$
are uniformly bounded. The desired result follows because $S[u_t]=S[w_t]$.

A standard application of Brower's fixed point Theorem (see for instance \cite[Lecture 3, Lemma 2]{Alice})
implies the existence of a conformal diffeomorphism  $T$  such that, for $i=1,2,3,$
\begin{equation}\label{centered}
 \int_{S^2}x_i\exp(-4w_t\circ F_t+4\gamma_t)d\mu_0=0,\quad\mbox{where}\quad F_t^*(g_0)=\exp(2\gamma_t)g_0.
\end{equation}
Moreover, $\gamma_t$ has an upper bound given by
\begin{lemm}\label{pixies}
There is a universal constant $C_0$ so that
$$ \sup\lvert\gamma_t \rvert\leq C_0+ 2/3(\overline{r}_t-\underline{r}_t).$$
\end{lemm}
\begin{proof}
    Assume that the maximum of $\alpha_t$ is attained at the north pole $p$ in $S^2$. We know that
    $$\int_{\{x_3\geq 0\}}x_3 \exp(-4w_t\circ T+4\gamma_t)d\mu_0=\int_{\{x_3\leq 0\}}(-x_3)
    \exp(-4w_t\circ T+4\gamma_t)d\mu_0$$ and so
    \begin{multline}\label{what}
        1/2\int_{\{x_3\geq 1/2\}} \exp(-4\sup{w_t}+4\gamma_t)d\mu_0\\
         \leq \int_{\{x_3\leq 0 \}} \exp(-4\inf {w_t}+4\gamma_t)d\mu_0.
        \end{multline}
    In stereographic coordinates,  $\gamma_t$ can be written as
    $$\gamma_t(x)=\log(\lambda(1+\lvert x\rvert^2)/(\lambda^2+\lvert x\rvert^2))$$
    for some $\lambda \leq 1$. Note that $\sup\gamma_t=-\inf\gamma_t=-\log\lambda$. An explicit computation shows
    the existence of some $r_0$ so that, if
        $$f(r)\equiv 8/3\lambda^4(\lambda^4+\lambda^2+3r^4+3(\lambda^2+1)r^2+1)/(\lambda^2+r^2)^3, $$
    then
        $$ \int_{\{x_3\geq 1/2\}} \exp(4\gamma_t)d\mu_0=f(0)-f(r_0)$$
    and
        $$ \int_{\{x_3\leq 0\}} \exp(4\gamma_t)d\mu_0=f(1).$$
    Hence, there are universal constants $C_1, C_2,$ and $C_3$ so that, for all $\lambda\leq 1$,
    $$\int_{\{x_3\geq 1/2\}} \exp(4\gamma_t)d\mu_0\geq C_1 \lambda^{-2}-C_2\lambda^{4}$$
    and
    $$\int_{\{x_3\leq 0\}} \exp(4\gamma_t)d\mu_0\leq C_3\lambda^{4}.$$
    Therefore, \eqref{what} implies that, for some  universal constants $C_5$
    $$\lambda^{-6}\leq C_5\exp(4\sup{w_t}-4\inf w_t)+C_5.$$
    The result follows because $\sup{w_t}-\inf w_t=\overline{r}_t-\underline{r}_t$.
\end{proof}

Let
$$u_t\equiv w_t\circ T-\gamma_t+\log\delta_t/2,\quad\mbox{where}\quad\delta_t\equiv\fint_{S^2}\exp(-2w_t)d\mu_0.$$
The effect of $\delta_t$ is to ensure that $\exp(-2u_t)g_0$ has volume $4\pi$ and we observe that, combining
Theorem \ref{intrinsic} with Proposition \ref{estimates}, we obtain
$$\delta_t=1+O(\exp(-\underline{r}_t))+O(\exp(2\overline{r}_t-3\underline{r}_t)).$$
Because
$$S[u_t]=S[w_t]-\log\delta_t,$$ it
suffices to show that $S[u_t]$ is uniformly bounded.

Because
$$\Delta_0 \gamma_t=1-\exp(2\gamma_t),$$
we have \eqref{e2} that
 $$\Delta_0 u_t=\exp(-2u_t)-1+Q_t,$$
 where $$Q_t=\exp(2\gamma_t)P_t\circ T+\exp(-2u_t)(\delta_t-1).$$
 Integration by parts yields
 \begin{multline*}
    2\int_{S^2}\lvert\nabla_0 u_t\rvert^2d\mu_0-2\int_{S^2}u_td\mu_0\\= \int_{S^2}-
    2u_t\exp(-2u_t)d\mu_0-\int_{S^2}2u_t Q_td\mu_0.
 \end{multline*}
We argue that the last term on the right-hand side is bounded independently of $t$. Note that Lemma \ref{pt} and
Lemma \ref{pixies} imply that
\begin{align*}
    \int_{S^2}\lvert u_t Q_t \rvert d\mu_0 &\leq\int_{S^2}\lvert u_t\circ T^{-1} P_t \rvert d\mu_0+\lvert\delta_t-1\rvert
    \int_{S^2}\lvert u_t\circ T^{-1}\rvert \exp(-2w_t) d\mu_0\\
    & \leq (\lvert\overline{r}_t-\underline{r}_t\rvert+C_0) \left( O(-\underline{r}_t)+O(\exp(2\overline{r}_t
    -3\underline{r}_t))\right)
\end{align*}
and so, because we are assuming condition \eqref{condi}, we obtain
$$\lim_{t\to\infty}\int_{S^2}\lvert u_t Q_t \rvert d\mu_0=0.$$

Thus, using Jensen's inequality we obtain that, for all $t$ sufficiently large,
\begin{equation}\label{alice}
2\int_{S^2}\lvert\nabla_0 u_t\rvert^2d\mu_0-2\int_{S^2}u_td\mu_0\leq 4\pi\log\fint_{S^2}\exp(-4u_t)d\mu_0+1.
\end{equation}

\begin{lemm}
$$\int_{S^2}\lvert\nabla_0 u_t\rvert d\mu_0$$ is uniformly bounded independently of $t$.
\end{lemm}
\begin{proof}
    Denoting by $G(x,y)$ the Green's function for the Laplacian on $S^2$ we have
    $$u_t(x)=\int_{S^2}u_td\mu_0-\int_{S^2}\Delta_0 u_t(y) \,G(x,y) d\mu_0(y) $$
    and so, because $\nabla G(\cdot,y)$ is in $L^1$ for all $y$,
    $$ \int_{S^2}\lvert\nabla_0 u_t\rvert d\mu_0\leq \int_{S^2}\lvert\Delta_0 u_t\rvert d\mu_0\leq C$$
    for some constant $C$ independent of $t$.
\end{proof}
This lemma and identity \eqref{centered} allow us to use an improvement of Aubin's inequality \cite{Aubin}, due
to Alice Chang and Paul Yang \cite[Lemma 4.2 ]{AlicePaul}, which says that under these conditions

$$4\pi\log\fint_{S^2}\exp(-4u_t)d\mu_0\leq 2\int_{S^2}\lvert\nabla_0 u_t\rvert^2 d\mu_0-4\int_{S^2}u_td\mu_0+C$$
where $C$ is a constant independent of $t$. Combining this with \eqref{alice} we get that
$$\hat{u}_t\equiv\int_{S^2}u_t\,d\mu_0$$ is uniformly bounded from above. A uniform bound below follows trivially
from $\exp(-2u_t)g_0$ having volume $4\pi$. In order to bound the gradient term we integrate by parts again so
that, for all $t$ sufficiently large,
\begin{align*}
    2\int_{S^2}\lvert\nabla_0 u_t\rvert^2d\mu_0 & =2\int_{S^2}\exp(-2u_t)(\hat{u}_t-u_t)d\mu_0+\int_{S^2}2(\hat{u}_t-u_t)
    Q_td\mu_0\\
    & \leq 2\int_{S^2}(\exp(-2u_t)-2\pi)(\hat{u}_t-u_t)d\mu_0+1\\
    & \leq  2\pi\log\left( (2\pi)^{-1}\int_{S^2}(\exp(-2u_t)-2\pi)\exp(2(\hat{u}_t-u_t))d\mu_0\right)\\
    &\quad+1\\
    & \leq 2\pi\log\int_{S^2}\exp(-4u_t)d\mu_0+C\\
    & \leq \int_{S^2}\lvert\nabla_0 u_t\rvert^2 d\mu_0+C,
\end{align*}
 where $C$ is a constant independent of $t$.

\newpage

\bibliographystyle{amsbook}

\vspace{20mm}

\end{document}